\documentclass[10pt, amstex]{article}
\usepackage{amssymb,amsmath}
\usepackage{mathrsfs}
\newcommand{\bbox}{\hfill $\Box$}

\begin{document}
\newcommand{\bea}{\begin{eqnarray}}
\newcommand{\ena}{\end{eqnarray}}
\newcommand{\beas}{\begin{eqnarray*}}
\newcommand{\enas}{\end{eqnarray*}}
\newcommand{\beq}{\begin{equation}}
\newcommand{\enq}{\end{equation}}
\newcommand{\qmq}[1]{\quad \mbox{#1} \quad}
\newcommand{\qm}[1]{\quad \mbox{#1}}
\def\qed{\hfill \mbox{\rule{0.5em}{0.5em}}}
\newcommand{\From}{From}
\newcommand{\ignore}[1]{}
\newcommand{\ws}{{\widetilde \sigma}}
\newcommand{\btheta}{\mbox{\boldmath {$\theta$}}}
\newcommand{\transpose}{{\mbox{\scriptsize\sf T}}}
\newcommand{\B}{B}
\newcommand{\E}{E}
\newcommand{\V}{V}
\newcommand{\Bvert}{\left\vert\vphantom{\frac{1}{1}}\right.}
\newtheorem{theorem}{Theorem}[section]
\newtheorem{corollary}{Corollary}[section]
\newtheorem{conjecture}{Conjecture}[section]
\newtheorem{proposition}{Proposition}[section]
\newtheorem{lemma}{Lemma}[section]
\newtheorem{definition}{Definition}[section]
\newtheorem{example}{Example}[section]
\newtheorem{remark}{Remark}[section]
\newtheorem{case}{Case}[section]
\newtheorem{condition}{Condition}[section]
\newcommand{\pf}{\noindent {\bf Proof:} }
\newcommand{\nn}{\nonumber}
\newcommand{\proof}{{\it Proof.}\ }

\title{{\bf\Large A Berry Esseen Theorem for the Lightbulb Process}}
\author{Larry Goldstein and Haimeng Zhang\\University of Southern California and Mississippi State University}
\footnotetext{AMS 2010 subject classifications: Primary:
62E17\ignore{Approximations to distributions (nonasymptotic)},
60C05\ignore{Combinatorial probability} Secondary: 60B15\ignore{Probability
measures on groups, Fourier transforms, factorization},
62P10\ignore{Applications to biology and medical sciences}.}
\footnotetext{Key words and phrases: Normal approximation, toggle
switch, Stein's method, size biasing.} \maketitle

\date{}

\begin{abstract}
In the so called lightbulb process, on days $r=1,\ldots,n$, out of $n$ lightbulbs, all initially off, exactly
$r$ bulbs, selected uniformly and independent of the past, have their status changed from off to on, or vice versa. With $X$ the number of bulbs on at the terminal time $n$, an even integer, and $\mu=n/2, \sigma^2=\mbox{Var}(X)$, we have
$$
\sup_{z \in \mathbb{R}} \left|P\left( \frac{X-\mu}{\sigma} \le z \right)-P(Z \le z)\right| \le \frac{n}{2\sigma^2}
\overline{\Delta}_0 + 1.64 \frac{n}{\sigma^3}+ \frac{2}{\sigma}
$$
where $Z$ is a standard normal random variable, and
$$
\overline{\Delta}_0 = \frac{1}{2\sqrt{n}} + \frac{1}{2n} + \frac{1}{3}e^{-n/2} \qmq {for $n \ge 6$,}
$$
yielding a bound of order $O(n^{-1/2})$ as $n \rightarrow \infty$.
A similar, though slightly larger bound holds for $n$ odd. The results are shown using a version of Stein's method for bounded, monotone size bias couplings. The argument for even $n$ depends on the construction of a variable $X^s$ on the same space as $X$ that has the $X$-size bias distribution, that is, that satisfies
\beas
E [X g(X)] =\mu E[g(X^s)] \quad \mbox{for all bounded continuous $g$,}
\enas
and for which there exists a $B \ge 0$, in this case $B=2$, such that $X \le X^s \le X+B$ almost surely. The argument for $n$ odd is similar to that for $n$ even, but one first couples $X$ closely to $V$, a symmetrized version of $X$, for which a size bias coupling of $V$ to $V^s$ can proceed as in the even case. In both the even and odd cases, the crucial calculation of the variance of a conditional expectation requires detailed information on the spectral decomposition of the lightbulb chain.
\end{abstract}

\section{Introduction}
The problem we consider here arises from a study in the pharmaceutical
industry on the effects of dermal patches designed to activate targeted receptors.
An active receptor will become inactive, and an inactive one active,
if it receives a dose of medicine released from the dermal patch. Let
the number of receptors, all initially inactive, be denoted by $n$.
On each day of the study, some number of randomly selected receptors will each receive one dose
of medicine, changing their statuses between the inactive and active states.
We adopt the following, somewhat more colorful, though equivalent, `lightbulb process' formulation from [\ref{rrz}].
Consider $n$ toggle switches, each connected to a lightbulb, all of which are initially off. Pressing the toggle
switch connected to a bulb changes its status from off to on and vice versa.
The problem of determining the properties of $X$, the number of light bulbs on at the end of day $n$, was first considered in [\ref{rrz}] for the case where on each day $r=1,\ldots,n,$ exactly $r$ of the $n$ switches are randomly pressed.

More generally, consider the lightbulb process on $n$ bulbs with some number $k$ of stages, where
$s_r \in \{0,\ldots,n\}$ lightbulbs are toggled in stage $r$, for
$r=1,\ldots,k$; we refer to the vector ${\bf s}=(s_1,\ldots,s_k)$ recording the
number of bulbs affected on each study day as the `switch pattern.'
In order
to consider quantities that depend on some subset of size $b$ of the $n$ bulbs, we define
\bea
\label{def-lam} \lambda_{n,b,s}= \sum_{t=0}^b  {b \choose t}
(-2)^t\frac{(s)_t}{(n)_t} \qmq{and} \lambda_{n,b,{\bf
s}}=\prod_{r=1}^k \lambda_{n,b,s_r},
\ena
where $(n)_k=n(n-1)\cdots(n-k+1)$ denotes the falling factorial, and the
empty product is 1.
Generalizing the results in [\ref{rrz}], writing $X_{\bf s}$
for the number of bulbs on at the terminal time when applying the
switch pattern ${\bf s}=(s_1,\ldots,s_n)$, the martingale method in
Proposition 4 of [\ref{zl}] shows that if the process is initialized
with all bulbs off, then
\bea \label{genmeanvar} EX_{\bf s} =
\frac{n}{2}\left( 1- \lambda_{n,1,{\bf s}} \right) \qmq{and}
\mbox{Var}(X_{\bf s})=\frac{n}{4}(1-\lambda_{n,2,{\bf s}})+
\frac{n^2}{4}(\lambda_{n,2,{\bf s}}-\lambda_{n,1,{\bf s}}^{\,2}),
\ena
where, from (\ref{def-lam}),
\bea
\label{lambdan12s}
\lambda_{n,1,s}=
1-\frac{2s}{n} \qmq{and} \lambda_{n,2,s}=1 -
\frac{4s}{n}+\frac{4s(s-1)}{n(n-1)} \quad \mbox{for $s=1,\ldots,n$.}
\ena

Letting ${\bf n}=(1,\ldots,n)$, we call the standard lightbulb
process the one where ${\bf s}={\bf n}$, and in this case
we will write $X$ short for
$X_{\bf n}$. In particular, (\ref{genmeanvar}) with ${\bf s}={\bf n}$ recovers the mean
$\mu=EX$ and variance of $\sigma^2=\mbox{Var}(X)$ as computed in
[\ref{rrz}]. Other results in [\ref{rrz}]
include recursions for determining the exact finite sample
distribution of $X$. Though computational approximations to the distribution of
$X$, including by the normal, were also considered in [\ref{rrz}],
the quality of such approximations, and the asymptotic normality of
$X$ was left open.

Theorem \ref{thm:main} settles the matter of the asymptotic
distribution of $X$ by providing a bound to the normal which holds
for all finite $n$, and which tends to zero at the rate $n^{-1/2}$
as $n$ tends to infinity.
We consider the case of even and odd $n$
separately. In the even case we directly couple the variable $X$
to one having the $X$-size bias distribution, as described later on in this section. In the even case (\ref{lambdan12s}) yields $\lambda_{n,1,n/2}=0$, and therefore $\lambda_{n,1,{\bf n}}=0$, hence from (\ref{genmeanvar}) we obtain $EX=n/2$, and also that $\sigma^2=\mbox{Var}(X)$ is given by
\bea \label{even.variance}
\sigma^2 = \frac{n}{4}(1-\lambda_{n,2,{\bf n}})+
\frac{n^2}{4}\lambda_{n,2,{\bf n}}.
\ena

To state our result for $n$ odd, for ${\bf s}=(s_1,\ldots,s_k)$ let
\bea \label{def:overlinelambda}
{\overline \lambda}_{n,b,r} = \left\{
\begin{array}{cc}
\frac{1}{2}\left(\lambda_{n,b,m} + \lambda_{n,b,m+1} \right) & r \in \{m,m+1\}\\
\lambda_{n,b,r} & \mbox{otherwise}
\end{array}
\right.
\qmq{and}
{\overline \lambda}_{n,b,{\bf s}}=\prod_{j=1}^k {\overline \lambda}_{n,b,s_j},
\ena
that is, ${\overline \lambda}_{n,b,{\bf s}}$ is obtained from $\lambda_{n,b,{\bf s}}$ by
replacing any occurrences of $\lambda_{n,b,m}$ and $\lambda_{n,b,m+1}$ in the product (\ref{def-lam}) by their average.
In the odd case, we proceed by first coupling $X$ to a more symmetric
random variable $V$ with mean and variance given respectively by
\bea \label{defmuV}
\mu_V=\frac{n}{2}\qmq{and}
\sigma_V^2
=\frac{n}{4}\left(1-{\overline \lambda}_{n,2,{\bf n}}\right)+\frac{n^2}{4}{\overline \lambda}_{n,2,{\bf n}}. \ena
Then, with $V$ in hand, we couple $V$ to a variable with the $V$-size bias distribution, and proceed as in the even case.
In Theorem \ref{thm:main}, and the remainder of the paper, $Z$ denotes a standard
normal random variable.
\begin{theorem}
\label{thm:main} Let $X$ be the number of bulbs on at the terminal
time $n$ in the standard lightbulb process.
Then for all $n$ even, with $\sigma^2$ as given in (\ref{even.variance}),
\bea \label{thm:evb}
\sup_{z \in \mathbb{R}} \left|P\left( \frac{X-n/2}{\sigma} \le z
\right)-P(Z \le z)\right| \le \frac{n}{2\sigma^2}{\overline \Delta}_0 + 1.64
\frac{n}{\sigma^3}+ \frac{2}{\sigma} \qmq
{for all $n \ge 6$,}
\ena
where
\beas
\overline{\Delta}_0 = \frac{1}{2\sqrt{n}} + \frac{1}{2n} + \frac{1}{3}e^{-n/2},
\enas
and for all $n=2m+1$ odd,
\beas
\sup_{z \in \mathbb{R}} \left|P\left( \frac{X-n/2}{\sigma_V}\le z
\right)-P(Z \le z)\right| \le \frac{n}{2\sigma_V^2}{\overline \Delta}_1 +  1.64
\frac{n}{\sigma_V^3}+ \frac{2}{\sigma_V}\left(1 + \frac{1}{\sqrt{2
\pi}}\right) \qmq {for all $n \ge 7$,}
\enas
where $\sigma_V^2$ is given in (\ref{defmuV}) and
\bea \label{upDelta1}
{\overline \Delta}_1 = \frac{1}{\sqrt{n}} +\frac{1}{2 \sqrt{2}}e^{-n/4}.
\ena
\end{theorem}
In the even case, as $\lambda_{n,2,{\bf n}}$ decays exponentially fast to zero, the variance $\sigma^2$ is of order $n$ and the bound (\ref{thm:evb}), therefore, of order $1/\sqrt{n}$; analogous remarks hold for the case where $n$ is odd.

We now more formally describe the lightbulb process on $n$ bulbs with $k$ stages. With $n \in
\mathbb{N}$ fixed and ${\bf s}=(s_1,\ldots,s_k)$ with $s_r \in
\{0,\ldots,n\}$ for $r=1,\ldots,k$, we will let ${\bf X}_{\bf
s}=\{X_{rj}: r=0,1,\ldots,k ,j=1,\ldots,n\}$ denote a collection of
Bernoulli variables. The initial state of the bulbs is given deterministically by $\{X_{0j},j=1,\ldots,n\}$, which will be taken to be state zero, that is, all bulbs off, unless specifically stated otherwise; non zero initial conditions are considered only in the Appendix. For $r \in \{1,\ldots,k\}$ the components of the `switch variables' ${\bf X}_{\bf s}$
have the interpretation that \beas X_{rj} &=& \left\{
\begin{array}{cc}
1 & \mbox{ if the status of bulb $j$ is changed at stage $r$,} \\
0 &\mbox{ otherwise.}
\end{array}
\right.
\enas

At stage $r$, $s_r$ of the $n$ bulbs are chosen uniformly to have their status changed, and the stages are independent of each other. Hence, with ${\bf e}=\{e_{rj}\}_{1 \le r \le k, 1 \le j \le n}$ an array of $\{0,1\}$ valued variables, the distribution of ${\bf X}_{\bf s}$ is given by
\bea \label{Xik-distribution}
P({\bf X}_{\bf s}={\bf e})
=\left\{
\begin{array}{cl}
\prod_{r=1}^k {n \choose s_r}^{-1} & \mbox{if $\sum_{j=1}^n e_{rj} = s_r, r=1,\ldots,k$}  \\
0 & \mbox{otherwise.}
\end{array}
\right.
\ena
Clearly, the vectors of stage $r$ switch variables $(X_{r1},\ldots,X_{rn})$ are exchangeable and the marginal distribution of the components $X_{rj}$ are Bernoulli with success probability $s_r/n$.
In general, for $j = 1, \ldots, n$,
the variables
\bea
\label{def-Yk}
X_j = \left(\sum_{r=0}^k X_{rj}\right) \mbox{ mod }2 \qmq{and} X_{\bf s}=\sum_{j=1}^n X_j
\ena
are the indicator that
bulb $j$ is on at the terminal time $k$, and the total number of bulbs on at that time, respectively.
For the standard lightbulb process we will write ${\bf X}$ and $X$ for ${\bf X}_{\bf n}$ and $X_{\bf n}$ respectively.

The lightbulb process, where the individual states of the $n$ bulbs evolve according
to the same marginal Markov chain, is a special case of a class of multivariate
chains studied in [\ref{zl}], known as
composition Markov chains of multinomial type. As shown in [\ref{zl}],
such chains admit explicit full spectral
decompositions, and in particular, the transition matrices for the stages of the lightbulb
process can be simultaneously diagonalized by a Hadamard matrix. These properties
were put to use in [\ref{rrz}] for the calculation of the
moments needed to compute the mean and variance of $X$. Here we put these same
properties to somewhat more arduous work, the calculation of moments of fourth order.

That no higher order moments are required for the derivation of a finite sample bound holding
for all $n$ is one distinct advantage of the technique we apply here, Stein's method for
the normal distribution, brought to life in the seminal
monograph [\ref{Stein86}]. By contrast,
the method of moments requires the calculation and appropriate convergence of moments of all orders,
and obtains only convergence in distribution. Stein's method for the normal is based on the
characterization of the normal distribution in [\ref{Steinchar}], which states that $Z$ is a
standard normal variable if and only if
\bea
\label{char:stein}
E[Z g(Z)]=E[g'(Z)]
\ena
for all absolutely continuous functions $g$ for which these expectations exist. The idea behind Stein's
method is that if a mean zero, variance one random variable $W$ is close in distribution to $Z$, then $W$
will satisfy (\ref{char:stein}) approximately. Hence, to gauge the proximity of $W$ to $Z$ for a given
test function $h$, one can evaluate the difference $Eh(W)-Nh$, where $Nh=Eh(Z)$, by solving the Stein equation
\beas
f'(w)-wf(w)=h(w)-Nh
\enas
for $f$ and evaluating $E[f'(W)-Wf(W)]$.
A priori it may appear that an evaluation of $E[f'(W)-Wf(W)]$ would
be more difficult than that for $Eh(W)-Nh$. However, the former form may be
handled through couplings.

Here we consider size bias
couplings to evaluate $E[f'(W)-Wf(W)]$. Given a nonnegative random variable $Y$
with positive finite mean $\mu=EY$, we say $Y^s$ has the $Y$-size
bias distribution if $P[Y^s \in dy] = (y/\mu) P[Y \in dy]$, or
more formally, if
\bea \E [Y g(Y)]=\mu \E [g(Y^s)] \quad \mbox{for all
bounded continuous functions $g$.} \label{formalsb}
\ena

The use of size bias couplings in Stein's method was introduced in
[\ref{BRS}], where it was applied to derive bounds of order
$\sigma^{-1/2}$ for the normal approximation to the number of local
maxima $Y$ of a random function on a graph, where
$\sigma^2=\mbox{Var}(Y)$. In [\ref{GR}] the method was extended to
multivariate normal approximations, and the rate was improved to
$\sigma^{-1}$, for the expectation of smooth functions of a vector
${\bf Y}$ recording the number of edges with certain fixed degrees
in a random graph. In [\ref{Gold}] the method was used to give
bounds in the Kolmogorov distance of order $\sigma^{-1}$ for various
functions on graphs and permutations, and in [\ref{GolPen}] for two
problems in the theory of coverage processes, with bounds of this
same order. A more complete treatment of Stein's method and its applications can
be found in [\ref{Chen}].

Here we prove and apply Theorem \ref{L-BD:SDKsConc} for bounded, monotone size
bias couplings, which
requires that the random variable $Y$ of interest, and a random variable $Y^s$, having the $Y$-size bias distribution, be constructed on a common space such that for some nonnegative constant $B$,
\beas
Y \le Y^s \le Y+B
\enas
with probability one.
Loosely speaking, Theorem \ref{L-BD:SDKsConc} says that
 given any such coupling of
$Y$ and $Y^s$ on a common space, an upper bound on the Kolmogorov
distance between the
distribution of $Y$ and the normal can computed in terms of $EY,\mbox{Var}(Y),B$ and
the quantity
\bea \label{def-Delta-int}
\Delta = \sqrt{\mbox{Var}(E(Y^s-Y|Y))}.
\ena
Theorem \ref{L-BD:SDKsConc} is based on a concentration
type inequality provided in Lemma \ref{L-BD:ConcSB}

For the standard lightbulb process, a size bias coupling of $X$ to $X^s$ is achieved in the even case by the construction, for each $i=1,\ldots,n$, of a collection
${\bf X}^i$ from the given ${\bf X}$ as follows. Recalling (\ref{def-Yk}), where
for ${\bf s}={\bf n}$ we have $k=n$, if $X_i=1$, that is, if bulb $i$ is on
at the terminal time, we set ${\bf X}^i={\bf X}$. Otherwise,
let $J$ be uniformly chosen from all $j$ for which $X_{n/2,j}=1-X_{n/2,i}$ and
let ${\bf X}^i$ be the same as ${\bf X}$ but that the values of
$X_{n/2,i}$ and $X_{n/2,J}$ are interchanged. Let $X^i$ be the number of bulbs
on at the terminal time when applying the switch variables ${\bf X}^i$.
Then, with $I$ uniformly chosen from $1,\ldots,n$, the
variable $X^s=X^I$ has the $X$ size bias distribution, essentially due to the fact,
shown in Lemma \ref{lem:even-given=1}, that
\beas
{\cal L}({\bf X}^i)={\cal L}({\bf X}|X_i=1).
\enas

Due to the parity issue, to handle the odd case when $n=2m+1$, we
first construct a coupling of $X$ to a more symmetric variable $V$.
The variable $V$ is constructed by randomizing stages $m$ and $m+1$ in the switch variables that
yield $X$. In particular at stage $m$ we add an additional switch with probability $1/2$ and, independently, at stage $m+1$ we remove an existing switch with probability $1/2$.
A size bias coupling of $V$ to $V^s$ can be
achieved as in the even case, thus yielding a bound to the normal
for $X$. We remark that the size biased couplings developed here are
used in [\ref{GhGo}] to show the distribution of $X$, in both the even and odd cases,
obey concentration of measure inequalities.

In Section \ref{sec:thm} we present Theorem \ref{L-BD:SDKsConc},
which gives a bound to the normal when a bounded, monotone size
biased coupling can be constructed for a given $X$. Our coupling
construction and the proof of the bound for the even case of the
lightbulb process is given in Section \ref{sec:ccn}.
Symmetrization, that is, the construction of $V$ from
$X$, coupling constructions for $V$, and the proof of the bound for
the odd case are given in Section \ref{sec:apv}. Calculations
of the bounds on the variance $\Delta$ in (\ref{def-Delta-int})
require estimates on $\lambda_{n,b,{\bf s}}$ in (\ref{def-lam}).
These estimates, based on the work of [\ref{zl}], yield the
spectral decomposition of the underlying transition matrices of the chain
and are given in Section \ref{sec:spe}. Some of the more
detailed calculations have been relegated to the Appendix.

\section{Bounded Monotone Couplings}
\label{sec:thm}
Theorem \ref{L-BD:SDKsConc} for bounded monotone size bias couplings depends
on the following lemma, which is in some sense the size bias version of Lemma 2.1 of [\ref{ShSu}]. With
$Y$ having mean $\mu$ and variance $\sigma^2$, both finite and positive, with some slight abuse of notation in
the definition of $W^s$, we set
\bea \label{def-WWs}
W=\frac{Y-\mu}{\sigma} \quad \mbox{and} \quad W^s=\frac{Y^s-\mu}{\sigma}.
\ena

\begin{lemma}
\label{L-BD:ConcSB}
Let $Y$ be a nonnegative random variable with mean $\mu$ and variance $\sigma^2$, both finite and positive,
and let $Y^s$ be
given on the same space as $Y$, having the $Y$-size bias distribution, and
satisfying $Y^s \ge Y$ with probability one. Then with $W$ and $W^s$ given in (\ref{def-WWs}),
for any $z \in \mathbb{R}$ and $a>0$,
\beas
\frac{\mu}{\sigma}E(W^s-W){\bf 1}_{\{W^s-W \le a\}}{\bf 1}_{\{z \le W \le z+a \}} \le a.
\enas
\end{lemma}

\noindent \proof For fixed $z \in \mathbb{R}$ let
\beas
f(w) = \left\{
\begin{array}{cl}
-a & w \le z\\
w-z-a & z < w \le z+2a\\
a & w > z+2a.
\end{array}
\right.
\enas
Then, using $|f(w)| \le a$ for all $w \in \mathbb{R}$, $\mbox{Var}(W)=1$ and the Cauchy Schwarz inequality to obtain the
first inequality, followed by definition (\ref{def-WWs}) and the size bias relation (\ref{formalsb}), we have
\beas
a &\ge& E(Wf(W))\\
&=&  \frac{1}{\sigma}E(Y-\mu)f\left(\frac{Y-\mu}{\sigma}\right)\\
&=& \frac{\mu}{\sigma} E\left(f(W^s)-f(W) \right) \\
&=& \frac{\mu}{\sigma} E\int_0^{W^s-W} f'(W+t)dt \\
&\ge& \frac{\mu}{\sigma} E\int_0^{W^s-W}{\bf 1}_{\{0 \le t \le a\}} {\bf 1}_{\{z \le W \le z+a\}} f'(W+t)dt,
\enas
where in the final inequality we have used that $W^s \ge W$ and $f'(w) \ge 0$ for all $w \in \mathbb{R}$.
Noting that $f'(W+t)={\bf 1}_{\{z \le W +t \le z+2a\}}$, and that $0 \le t \le a$ and $z \le W \le z+a$
imply $z \le W+t \le z+2a$, we have
$$
{\bf 1}_{\{0 \le t \le a\}} {\bf 1}_{\{z \le W \le z+a\}}f'(W+t)={\bf 1}_{\{0 \le t \le a\}}{\bf 1}_{\{z \le W \le z+a\}},
$$
and therefore obtain
\beas
a &\ge&\frac{\mu}{\sigma} E \int_0^{W^s-W} {\bf 1}_{\{0 \le t \le a\}}{\bf 1}_{\{z \le W \le z+a\}}dt \\
&=& \frac{\mu}{\sigma} E\left(\min(a,W^s-W) {\bf 1}_{\{z \le W \le z+a\}}\right)\\
&\ge& \frac{\mu}{\sigma} E(W^s-W) {\bf 1}_{\{W^s-W \le a\}}{\bf 1}_{\{z \le W \le z+a\}},
\enas
as claimed.
\bbox

\begin{theorem}
\label{L-BD:SDKsConc}
Let $Y$ be a nonnegative random variable with mean $\mu$ and  variance $\sigma^2$, both finite and positive, and let $Y^s$ be given on the same space as $Y$, with the $Y$-size bias distribution, satisfying $Y \le Y^s \le Y+B$ with probability
one, for some positive constant $B$. Then with $W$ and $W^s$ given in (\ref{def-WWs}),
we have
\beas
\sup_{z \in \mathbb{R}}|P(W \le z) - P(Z\le z)|
\le \frac{\mu}{\sigma^2}\Delta + 0.82 \frac{\delta^2 \mu}{\sigma}+ \delta,
\enas
where
\bea
\label{def-Delta}
\Delta = \sqrt{\mbox{Var}(E (Y^s-Y|Y))} \quad \mbox{and} \quad \delta=B/\sigma .
\ena
\end{theorem}

\noindent \proof For $z \in \mathbb{R}$ arbitrary, let $h(w)={\bf 1}_{\{w \le z\}}$ and
let $f(w)$ be the unique bounded solution to the Stein equation
\bea
\label{stein2}
f'(w)-wf(w)=h(w)-Nh,
\ena
where $Nh=Eh(Z)$. Substituting $W$ into (\ref{stein2}) and using definition (\ref{def-WWs}) and
the size bias relation (\ref{formalsb}) yields
\bea
\nn
\lefteqn{E\left( h(W)-Nh \right)} \\
\nn &=& E\left(f'(W)-Wf(W) \right)\\
\nn &=& E\left(f'(W)- \frac{\mu}{\sigma}(f(W^s)-f(W))\right)\\
\label{condexptermpint}
&=&E\left(f'(W)(1-\frac{\mu}{\sigma}(W^s-W))-\frac{\mu}{\sigma}\int_0^{W^s-W}(f'(W+t)-f'(W))dt
\right).
\ena
As compiled in [\ref{Chen}], we have the following bounds on the solution $f$ from Lemma 2 in Chapter II of [\ref{Stein86}],
\bea
\label{bd:st86}
0<f(w)<\sqrt{2\pi}/4 \qmq{and} |f'(w)| \le 1,
\ena
and also, as previously noted in [\ref{ShSu}], as a consequence of (\ref{bd:st86}) and the mean value theorem,
we obtain
\bea
\label{bnd:ChSu}
|(w+t)f(w+t)-wf(w)| \le (|w|+\sqrt{2\pi}/4)|t|.
\ena

Noting that $EY^s=EY^2/\mu$ by (\ref{formalsb}), we find
$$
\frac{\mu}{\sigma}E(W^s-W)=\frac{\mu}{\sigma^2}\left( \frac{EY^2}{\mu}-\mu
\right)=1.
$$
Therefore, taking expectation by conditioning and then applying
(\ref{bd:st86}) and the Cauchy Schwarz inequality, we bound the first term in (\ref{condexptermpint}) as
\beas
\left| E  \left\{ f'(W) E  \left( 1 - \frac{ \mu }{\sigma}
(W^s-W)|W\right)
\right\}\right| \le \frac{ \mu }{\sigma} \sqrt{\mbox{Var}(E
(W^s-W|W))}= \frac{ \mu }{\sigma^2}\Delta.
\enas

To bound the remaining term of (\ref{condexptermpint}), using (\ref{stein2}), we have
\bea
\nn
\lefteqn{\frac{\mu}{\sigma}\int_0^{W^s-W}\left(f'(W+t)-f'(W)\right)dt}\\
&=& \frac{\mu}{\sigma}
\int_0^{W^s-W}[(W+t)f(W+t)-Wf(W)]dt
\label{L-BD:thmSBKJs}
+\frac{\mu}{\sigma}\int_0^{W^s-W}({\bf 1}_{\{W+t \le
z\}}-{\bf 1}_{\{W \le z\}})dt.
\ena
Applying (\ref{bnd:ChSu}) to the first term in (\ref{L-BD:thmSBKJs}), and using
$0 \le W^s-W \le \delta$ and $EW^2=1$, shows
that the absolute value of the expectation of this term is bounded by
\beas
\frac{\mu}{\sigma} E \left(
\int_0^{W^s-W}\left(|W|+\frac{\sqrt{2\pi}}{4}\right)t dt \right)
=  \frac{\mu}{2\sigma} E \left((W^s-W)^2 (|W|+\frac{\sqrt{2\pi}}{4})\right)
\le \frac{\mu}{2\sigma}\delta^2 (1+\frac{\sqrt{2\pi}}{4})
\le 0.82 \frac{\delta^2 \mu}{\sigma}.
\enas
Taking the expectation of the absolute value of the second term in
(\ref{L-BD:thmSBKJs}), we obtain
\beas
\frac{\mu}{\sigma}E \left| \int_0^{W^s-W}({\bf 1}_{\{W+t \le
z\}}-{\bf 1}_{\{W \le z\}})dt \right|
&=&\frac{\mu}{\sigma}E \left( \int_0^{W^s-W}{\bf 1}_{\{z-t<W\le
z\}}dt \right) \\
&\le&\frac{\mu}{\sigma}E (W^s-W){\bf 1}_{\{z-\delta<W\le z\}}
\enas
again using that $0 \le W^s-W \le \delta$ with probability 1.
Lemma \ref{L-BD:ConcSB} with $a=\delta$ and $z$ replaced by $z-\delta$ shows this
term can be no more than $\delta$. Since
$z \in \mathbb{R}$ was arbitrary
the proof is complete.
\bbox

\section{Normal approximation of $X$}
The next lemma shows that the size bias distribution of a sum may be achieved by taking certain mixtures.
The result is a special case of Lemma 2.1 of [\ref{GR}],
but we  give a short direct proof to make the paper more self-contained.
\begin{lemma}
\label{sblem}
Suppose $X$ is a sum of nontrivial exchangeable
Bernoulli variables $X_1,\ldots,X_n$, and that for $i \in \{1,\ldots,n\}$
the variables
$X^i_1,\ldots,X^i_n$ have joint distribution
$$
{\cal L}(X^i_1,\ldots,X^i_n)={\cal L}(X_1,\ldots,X_n|X_i=1).
$$
Then
$$
X^i=\sum_{j=1}^n X^i_j
$$
has the $X$-size bias distribution $X^s$, as does the mixture $X^I$ when $I$ is a random index with values in $\{1,\ldots,n\}$, independent of all other variables.
\end{lemma}

\noindent \proof For $i \in \{1,\ldots,n\}$, we first need to show that $X^i$ satisfies
(\ref{formalsb}), that is, that $E[X]E [g(X^i)] = E[X g(X)]$ holds for
a given bounded continuous $g$. Now, for such $g$ \beas E[X g(X)] =
\sum_{j=1}^n E[X_j g(X)]= \sum_{j=1}^n P[X_j=1]E[g(X)| X_j = 1].
\enas As exchangeability implies that $E[g(X)| X_j = 1]=E[g(X)| X_i = 1]$ for all $j=1,\ldots,n$, we have
\beas
E[X g(X)] =\left( \sum_{j=1}^n P[X_j=1]
\right)E[g(X)| X_i = 1] = E[X ] E[g(X^i)], \enas proving the
first claim.  The second claim now follows from
\beas
Eg(X^I)=\sum_{i=1}^n E[g(X^I),I=i] &=& \sum_{i=1}^n E[g(X^I)|I=i]P(I=i)= \sum_{i=1}^n E g(X^i)P(I=i)\\
&=&\sum_{i=1}^n E g(X^s)P(I=i)=Eg(X^s) \sum_{i=1}^n P(I=i) = Eg(X^s).
\enas
\bbox

\subsection{Even case}
\label{sec:ccn}

In this subsection we provide the proof of Theorem \ref{thm:main}
for even $n$. We begin by describing a coupling of $X$, the total number
of bulbs on at the terminal time $n$ in the standard lightbulb process, to a variable $X^s$ with the
$X$-size bias distribution. Throughout, we let ${\cal U}(S)$ denote the uniform
distribution over a finite set $S$.

\begin{theorem}
\label{thm:evncop}
With $n \in \mathbb{N}$ even, let the collection of switch variables ${\bf X}=\{X_{rj}: r,j=1,\ldots,n\}$ and $X$ satisfy (\ref{Xik-distribution}) and (\ref{def-Yk}), respectively, with ${\bf s}={\bf n}$.
For every $i=1,\ldots,n$ let ${\bf X}^i$ be given from ${\bf X}$ as follows. If $X_i=1$ then ${\bf X}^i={\bf X}$. Otherwise, with ${\cal L}(J^i|{\bf X}) = {\cal U}\{j: X_{n/2,j}=1-X_{n/2,i}\}$, let ${\bf X}^i=\{X_{rj}^i:r,j=1,\ldots,n\}$ where
\beas
X_{rj}^i=\left\{
\begin{array}{cl}
X_{rj} & r \not = n/2\\
X_{n/2,j}& r=n/2, j \not \in \{i,J^i\}\\
X_{n/2,J^i} & r=n/2, j=i\\
X_{n/2,i}& r=n/2, j=J^i,
\end{array}
\right.
\enas
and let $X^i=\sum_{j=1}^n X_j^i$ where
$$
X_j^i = \left(\sum_{r=1}^n X_{rj}^i\right) \mbox{ mod }2.
$$
Then, with $I$ uniformly
chosen from $\{1,\ldots,n\}$ and independent of all other variables, the mixture $X^I=X^s$ has the $X$-size
bias distribution and satisfies
\bea \label{YYsbdmn}
X^s-X=2 {\bf 1}_{\{X_I=0,
X_{J^I}=0\}} \qmq{and} X \le X^s \le X+2.
\ena
\end{theorem}

We prove Theorem \ref{thm:evncop} making use of a preliminary lemma, and also of the fact that
\bea \label{PYi=1/2}
P(X_j=0)=P(X_j=1)=\frac{1}{2} \quad \mbox{for all $j=1,\ldots,n$.}
\ena
The equalities in (\ref{PYi=1/2}) follow from $EX=n/2$, itself implied by (\ref{genmeanvar}) and that
$\lambda_{1,n,{\bf n}}=0$, as noted earlier.

\begin{lemma}
\label{lem:even-given=1}
For all $i=1,\ldots,n$, the collection of random variables ${\bf X}^i$ constructed from ${\bf X}$ as specified in Theorem \ref{thm:evncop} satisfies
\beas
{\cal L}({\bf X}^i)={\cal L}({\bf X}|X_i=1).
\enas
\end{lemma}

\noindent \proof For a given $i \in \{1,\ldots,n\}$ let ${\bf e}=\{e_{rj}: r,j=1,\ldots,n\}$ with $e_{rj} \in \{0,1\}$
for $r,j=1,\ldots,n$. First note that since ${\bf X}^i={\bf X}$ when $X_i=1$ we have
\beas
P({\bf X}^i={\bf e})&=&P(X_i=1)P({\bf X}^i={\bf e}|X_i=1)+P(X_i=0)P({\bf X}^i={\bf e}|X_i=0)\\
&=&P(X_i=1)P({\bf X}={\bf e}|X_i=1)+P(X_i=0)P({\bf X}^i={\bf e}|X_i=0),
\enas
so the desired conclusion is equivalent to
\bea \label{Xj|0}
P({\bf X}^i={\bf e}|X_i=0)=P({\bf X}={\bf e}|X_i=1).
\ena

As the construction of ${\bf X}^i$ preserves the number of switches in each stage $r$
we may assume $\sum_j e_{rj}=r$ for all $r$, as otherwise both sides of $(\ref{Xj|0})$ are zero.
If $\sum_r e_{ri}=0\,\mbox{mod}\,2$ then the left hand side of (\ref{Xj|0}) is zero since
$X_i^i=1$ by construction; similarly,
the right hand side is zero as ${\bf X}={\bf e}$ implies $X_i=0$.
Hence we need only verify (\ref{Xj|0}) assuming
\bea
\label{sumeri=1}
\sum_{j=1}^n e_{rj}=r \qmq{for all $r=1,\ldots,n$, and} \sum_{r=1}^n e_{ri} = 1\,\mbox{mod}\,2.
\ena

Writing $J$ for $J^i$ for simplicity, and letting ${\bf e}^{i,j}$ denote the array ${\bf e}$ with
coordinates $e_{n/2,i}$ and $e_{n/2,j}$ interchanged, by (\ref{PYi=1/2}) we have
\beas
P({\bf X}^i={\bf e}|X_i=0)&=& 2P({\bf X}^i={\bf e},X_i=0)\\
                          &=& 2\sum_{j=1}^n P({\bf X}^i={\bf e},X_i=0,J=j)\\
                            &=& 2\sum_{j=1}^n P({\bf X}={\bf e}^{i,j},X_i=0,J=j).
\enas
Note that when $e_{n/2,i}=e_{n/2,j}$, or equivalently, $e_{n/2,i}^{i,j}=e_{n/2,j}^{i,j}$ then
\bea \label{eij:bothsides}
P({\bf X}={\bf e}^{i,j},X_i=0,J=j)= P({\bf X}={\bf e}^{i,j},J=j)
\ena
as both sides are zero, since on $J=j$ we have $X_{n/2,i} \not = X_{n/2,j}$. Otherwise, $e_{n/2,i}\not =e_{n/2,j}$,
and by the second equality in (\ref{sumeri=1}),
\beas
\sum_{r=1}^n e_{ri}^{i,j} = \sum_{r\not = n/2} e_{ri}+ e_{n/2,i}^{i,j} = \sum_{r\not = n/2} e_{ri}+ 1-e_{n/2,i}
= \sum_{r=1}^n e_{ri}+ 1 = 0, \quad \mbox{with equalities modulo 2,}
\enas
so (\ref{eij:bothsides}) holds again. Hence
\beas
\lefteqn{P({\bf X}^i={\bf e}|X_i=0)}\\
&=& 2 \sum_{j=1}^n P({\bf X}={\bf e}^{i,j},J=j) \\
&=& 2 \sum_{j=1}^n P(J=j|{\bf X}={\bf e}^{i,j})P({\bf X}={\bf e}^{i,j})\\
&=&2\prod_{s=1}^n {n \choose s}^{-1}  \sum_{j=1}^n P(J=j|{\bf X}={\bf e}^{i,j})
=2\prod_{s=1}^n {n \choose s}^{-1},
\enas
where in the second to last equality we have used that ${\bf e}^{i,j}$, as ${\bf e}$, satisfies the first equality of (\ref{sumeri=1}), and the distribution of ${\bf X}$ given by (\ref{Xik-distribution}),
and in the last equality that the sum of probabilities of the conditional distribution of $J$ given an ${\bf X}$
configuration that satisfies the first equality of (\ref{sumeri=1}) for $r=n/2$ must sum to one.
Now, again using the second equality in (\ref{sumeri=1}),
\beas
2\prod_{s=1}^n {n \choose s}^{-1}= 2P({\bf X}={\bf e})=2P({\bf X}={\bf e},X_i=1)=P({\bf X}={\bf e}|X_i=1),
\enas
proving (\ref{Xj|0}), and the lemma. \bbox

We now present the proof of Theorem \ref{thm:evncop}.

\noindent \proof  That $X^s$ has the $X$-size bias distribution follows from Lemmas \ref{sblem} and \ref{lem:even-given=1}.
To prove the first equality in (\ref{YYsbdmn}), note that if $X_I=1$ then ${\bf X}^I={\bf X}$, hence in this case $X^s=X$.
Otherwise $X_I=0$ and the collection ${\bf X}^I$ is constructed from ${\bf X}$ by interchanging the stage $n/2$, unequal, switch variables $X_{n/2,I}$ and $X_{n/2,J^I}$. If $X_{J^I}=1$ then after the interchange $X_I^I=1$ and $X_{J^I}^I=0$, yielding $X^s=X$. If $X_{J^I}=0$ then after the interchange $X_I^I=1$ and $X_{J^I}^I=1$, yielding $X^s=X+2$. The second claim in (\ref{YYsbdmn}) is an immediate consequence of the first.
\bbox


The following lemma shows that for the case at hand the variance of
the conditional expectation term (\ref{def-Delta}) in Theorem
\ref{L-BD:SDKsConc} may be expressed in terms of quantities of the
form
\bea
\label{def:galphabeta}
g_{\alpha,\beta,{\bf s}}^{(l)}=P(X_1=
\cdots =X_{\alpha+\beta}=0, X_{l1}=\cdots =X_{l\alpha}=0,
X_{l,\alpha+1}=\cdots = X_{l,\alpha+\beta}=1),
\ena
the probability
that when applying the switch pattern ${\bf s}=(s_1,\ldots,s_k)$ on
$n$ bulbs over $k$ stages, bulbs numbered $1$ though $\alpha+\beta$
terminate in the off position, and in some stage $l \in \{1,\ldots,k\}$,
bulbs $1$ through $\alpha$ receive switch variable $0$, and bulbs
numbered $\alpha+1$ through $\alpha+\beta$ receive switch variable
$1$. We keep the number of bulbs $n$ implicit in the notation, and also suppress ${\bf s}$ when ${\bf s}={\bf n}$,
writing more simply $g_{\alpha,\beta}^{(l)}$.

Using the spectral decomposition in Section \ref{sec:spe} to handle
the probabilities in (\ref{def:galphabeta}), we now provide an upper bound to the
term (\ref{def-Delta}) when applying Theorem \ref{L-BD:SDKsConc} for $n$ even.
For a given ${\bf s}=(s_1,\ldots,s_k)$ and $l \in \{1,\ldots,k\}$ let
\bea
\label{sdell}
{\bf s}_l = (s_1,\ldots,s_{l-1},s_{l+1},\ldots,s_k),
\ena
the vector ${\bf s}$ with its $l^{th}$ component deleted. For $l \not = j$ we similarly let ${\bf s}_{l,j}$ denote
${\bf s}$ with its $l^{th}$ and $j^{th}$ components deleted.

\begin{lemma}
\label{lem:Delta0}
Let $n$ be even and $X$ and $X^s$ given by Theorem \ref{thm:evncop}. Then for $n \ge 6$,
\beas
\Delta_0 \le \overline{\Delta}_0  \qmq {where}
\Delta_0 = \sqrt{\mbox{Var}(E(X^s-X|X))} \qmq{and} \overline{\Delta}_0= \frac{1}{2\sqrt{n}} + \frac{1}{2n} + \frac{1}{3}e^{-n/2}.
\enas
\end{lemma}

\noindent \proof We apply the construction of $X^s$, and the conclusions, of Theorem \ref{thm:evncop}. For notational simplicity let $J^I=J$, so in particular, from (\ref{YYsbdmn}) we have $X^s-X=2{\bf 1}_{\{X_I=0,X_J=0\}}$. Expanding the indicator over the possible values of $I$ and $J$, and then over the values of the switch variables $X_{n/2,i}$ and $X_{n/2,j}$ yields
\beas
{\bf 1}_{\{X_I=0, X_J=0\}}&=& \sum_{i,j=1}^n {\bf 1}_{\{X_i=0, X_j=0\}}{\bf 1}_{\{I=i,J=j\}}\\
&=& \sum_{i,j=1}^n  {\bf 1}_{\{X_i=0, X_j=0, X_{n/2,i}=0\}}{\bf 1}_{\{I=i,J=j\}}+\sum_{i,j=1}^n{\bf 1}_{\{X_i=0, X_j=0, X_{n/2,i}=1\}}{\bf 1}_{\{I=i,J=j\}}\\
&=& \sum_{i \not = j} {\bf 1}_{\{X_i=0, X_j=0, X_{n/2,i}=0,X_{n/2,j}=1\}}{\bf 1}_{\{I=i,J=j\}}+\sum_{i \not = j}{\bf 1}_{\{X_i=0, X_j=0, X_{n/2,i}=1,X_{n/2,j}=0\}}{\bf 1}_{\{I=i,J=j\}}\\
&=&  2\sum_{i \not = j} {\bf 1}_{\{X_i=0, X_j=0, X_{n/2,i}=0,X_{n/2,j}=1\}}{\bf 1}_{\{I=i,J=j\}},
\enas
where the second to last equality holds almost surely,
as the probability of the event $\{I=i,J=j\}$ is zero whenever $X_{n/2,i}$ and $X_{n/2,j}$ agree, and the
last inequality holds as the final expression is the sum of two terms which can be seen to be equal
by reversing the roles of $i$ and $j$.

To obtain a tractable bound on the required variance we apply the inequality
\bea \label{cond.F}
\mbox{Var}(E(X^s-X|X)) \le \mbox{Var}(E(X^s-X|{\cal F})),
\ena
which holds when
${\cal F}$ is any $\sigma$-algebra with respect to which $X$ is measurable (see [\ref{GR}], for example). Here we let ${\cal F}$ be the $\sigma$-algebra generated by ${\bf X}$, the collection of all switch variables. The first indicator in the final sum
above, ${\bf 1}_{\{X_i=0, X_j=0, X_{n/2,i}=0,X_{n/2,j}=1\}}$, is
measurable with respect to ${\cal F}$. For the second indicator,
conditioning on ${\cal F}$ yields \beas
E \left({\bf 1}_{\{I=i,J=j\}}|{\cal F} \right)=
P(I=i,J=j|{\cal F})=\frac{2}{n^2}{\bf
1}_{\{X_{n/2,i} \not = X_{n/2,j}\}}, \enas as for any $i$, chosen with probability $1/n$, there are
$n/2$ choices for $j$ satisfying the condition in the indicator. Hence, recalling from (\ref{YYsbdmn})
that $X^s-X=2 {\bf 1}_{\{X_I=0,X_J=0\}}$, we
have
\bea \label{EYsminY} E\left(X^s - X \Bvert {\cal F}\right)&=&
U_n  \qmq{where} U_n = \frac{4}{n^2} \sum_{i \not =j} {\bf
1}_{\{X_i=0,X_j=0 , X_{n/2,i}=0,X_{n/2,j}=1\}},
\ena
and $\Delta_0^2 \le \mbox{Var}(U_n)$ by (\ref{cond.F}).

Taking the
expectation of $U_n$ in (\ref{EYsminY}), by the exchangeability of
the $(n)_2$ terms in the sum,  applying Corollary \ref{cor:spectral.4} and recalling the notation defined in (\ref{sdell}), we obtain
\bea
\label{U-expectation}
EU_n= \frac{4}{n^2}(n)_2 g_{1,1}^{(n/2)}= \frac{1}{4}\left(1-\lambda_{n,2,{\bf n}_{n/2}}\right).
\ena
Squaring (\ref{EYsminY}) in order to obtain the second moment of $U_n$, we
obtain a sum over indices $i_1,i_2,j_1,j_2$ with $\{i_1,i_2\} \cap \{j_1,j_2\} = \emptyset$, so $|\{i_1,i_2,j_1,j_2\}| \in \{2,3,4\}$, and we may write
\bea \label{UissumVs}
\lefteqn{U_n^2= U_{n,2}^2+U_{n,3}^2+U_{n,4}^2 \qmq{where}}\\
\nn &&U_{n,p}^2=\frac{16}{n^4}\sum_{\stackrel{|\{i_i,i_2,j_1,j_2\}|=p}{\{i_1,i_2\} \cap \{j_1,j_2\} = \emptyset}}{\bf 1}_{\{X_{i_1}=0,X_{j_1}=0 , X_{n/2,i_1}=0,X_{n/2,j_1}=1\}}{\bf 1}_{\{X_{i_2}=0,X_{j_2}=0 , X_{n/2,i_2}=0,X_{n/2,j_2}=1\}}.
\ena

Beginning the calculation with the main term $U_{n,4}^2$, where all four indices are distinct, taking expectation
using exchangeability, Corollary \ref{cor:spectral.4} yields
\bea \label{Vn4}
EU_{n,4}^2 = \frac{16}{n^4}(n)_4 g_{2,2}^{(n/2)} = \left( \frac{n-2}{4 n}\right)^2 \left( 1- 2\lambda_{n,2,{\bf n}_{n/2}}+\lambda_{n,4,{\bf n}_{n/2}}\right).
\ena

With the inequalities over the summation in (\ref{UissumVs})
in force, the event
$|\{i_1,i_2,j_1,j_2\}|=3$ can only occur when
\beas
a)\,i_1 \not = i_2, j_1=j_2\qmq{or} b)\, i_1=i_2, j_1 \not = j_2.
\enas
Now continuing to use Corollary \ref{cor:spectral.4} without further notice, case a) leads to a contribution of
\beas
\frac{16}{n^4}(n)_3 g_{2,1}^{(n/2)}=\frac{n-2}{4n^2}\left( 1+\lambda_{n,1,{\bf n}_{n/2}}-\lambda_{n,2,{\bf n}_{n/2}}-\lambda_{n,3,{\bf n}_{n/2}}\right),
\enas
while in the same manner the contribution from case b) is
\beas
\frac{16}{n^4}(n)_3 g_{1,2}^{(n/2)}=\frac{n-2}{4n^2}\left( 1-\lambda_{n,1,{\bf n}_{n/2}}-\lambda_{n,2,{\bf n}_{n/2}}+\lambda_{n,3,{\bf n}_{n/2}}\right).
\enas
Totaling we find
\bea \label{Vn3}
EU_{n,3}^2=\frac{n-2}{2n^2}\left( 1-\lambda_{n,2,{\bf n}_{n/2}}\right).
\ena

With the inequalities over the summation in (\ref{UissumVs})
in force, the event $|\{i_1,i_2,j_1,j_2\}|=2$ can
only occur when $i_1=i_2$ and $j_1=j_2$. Hence,
\bea \label{Vn2}
EU_{n,2}^2 = \frac{16}{n^4}(n)_2 g_{1,1}^{(n/2)}= \frac{1}{n^2}\left(1 - \lambda_{n,2,{\bf n}_{n/2}} \right).
\ena

Summing (\ref{Vn4}), (\ref{Vn3}) and (\ref{Vn2}) we obtain
\beas
EU_n^2 &=&\left( \frac{n-2}{4 n}\right)^2 \left( 1- 2\lambda_{n,2,{\bf n}_{n/2}}+\lambda_{n,4,{\bf n}_{n/2}}\right)+
\frac{1}{2n}\left(1 - \lambda_{n,2,{\bf n}_{n/2}} \right).
\enas
Subtracting the square of the first moment, given in (\ref{U-expectation}), yields
\beas
\mbox{Var}(U_n) &=& \frac{1}{16}\left(1- \frac{2}{n}\right)^2 \left( 1- 2\lambda_{n,2,{\bf n}_{n/2}}+\lambda_{n,4,{\bf n}_{n/2}}\right)+
\frac{1}{2n}\left(1 - \lambda_{n,2,{\bf n}_{n/2}} \right) -\frac{1}{16}\left(1-\lambda_{n,2,{\bf n}_{n/2}}\right)^2\\
&=&\frac{1}{16}\left(\lambda_{n,4,{\bf n}_{n/2}}- \lambda_{n,2,{\bf n}_{n/2}}^2\right)+ \frac{1-n}{4n^2} \left( 1- 2\lambda_{n,2,{\bf n}_{n/2}}+\lambda_{n,4,{\bf n}_{n/2}}\right)+
\frac{1}{2n}\left(1 - \lambda_{n,2,{\bf n}_{n/2}} \right)\\
&=&\frac{1}{4n}+\frac{1}{4n^2}+ \frac{1}{16}\left(\lambda_{n,4,{\bf n}_{n/2}}
- \lambda_{n, 2, {\bf n}_{n/2}}^2\right) - \frac{1}{2n^2} \lambda_{n, 2, {\bf n}_{n/2}} +
\frac{1-n}{4n^2}\lambda_{n, 4, {\bf n}_{n/2}}.
\enas

Now applying Lemma \ref{lam24bds}
from Section \ref{sec:spe}, for $n \ge 6$ we obtain
\beas \mbox{Var}(U_n) &\le&
\frac{1}{4n}+\frac{1}{4n^2}+
\frac{1}{16}\left(\frac{1}{2}e^{-n}+e^{-2n}\right) + \frac{1}{2n^2}
e^{-n} + \frac{n-1}{8n^2}e^{-n}\\
&\le& \frac{1}{4n}+\frac{1}{4n^2}+ e^{-n}\left(
\frac{1}{16} + \frac{1}{n^2} + \frac{1}{8n}\right). \enas
The inequality $\sqrt{a+b+c} \le \sqrt{a}+\sqrt{b}+\sqrt{c}$, holding
for all nonnegative $a,b$ and $c$, now yields the claim of the lemma. \bbox

With all ingredients at hand, we may now prove the bound for even $n$.

\noindent {\em Proof of Theorem \ref{thm:main}, even case.}
The size
biased coupling given in Theorem \ref{thm:evncop} satisfies the
hypotheses of Theorem \ref{L-BD:SDKsConc} with $B=2$, by the second
inequality in (\ref{YYsbdmn}). Hence the result
for the even case follows by applying Theorem \ref{L-BD:SDKsConc}
with $\mu=n/2,\delta=2/\sigma$ and the bound $\overline{\Delta}_0$ on $\Delta_0$ given in Lemma \ref{lem:Delta0}.
\bbox

\subsection{Odd case}
\label{sec:apv}
Now we move to the case where $n=2m+1$ is odd.
Instead of directly forming a size biased coupling to $X$, we first couple $X$ closely to a more symmetrical random variable $V$ for which a coupling like the one in the even case may be applied. The variable $V$ is constructed by randomizing stages $m$ and $m+1$. In particular at stage $m$ we add an additional switch with probability $1/2$ and, independently at stage $m+1$ we remove an existing switch with probability $1/2$.

Formally, let ${\bf X}=\{X_{rj}: r,j=1,\ldots,n\}$ be a collection of switch variables with distribution given by (\ref{Xik-distribution}) with ${\bf s}={\bf n}$, and let $X=X_{\bf n}$ be given by (\ref{def-Yk}) with $k=n$.
Let
\beas
{\cal L}(B_m|{\bf X}) = {\cal U}\{j:X_{mj}=0\} \qmq{and} {\cal L}(B_{m+1}|{\bf X}) = {\cal U}\{j:X_{m+1,j}=1\},
\enas
and let $C_m$ and $C_{m+1}$ be symmetric Bernoulli variables,
independent of ${\bf X}$ and of $B_m$ and $B_{m+1}$. Now let a collection of switch variables ${\bf V} = \{V_{rj}, r, j = 1, \ldots, n\}$ be defined by
\bea
\label{def:VX}
V_{rj} = \left\{
\begin{array}{cc}
X_{rj}  & r \not \in \{m,m+1\}\\
X_{mj} & r=m, j \not = B_m\\
C_m    & r=m,j=B_m\\
X_{m+1,j} & r=m+1,j \not = B_{m+1}\\
C_{m+1} & r=m+1,j=B_{m+1},
\end{array}
\right.
\ena
and set
\bea \label{defVVk}
V=\sum_{j=1}^n V_j \quad \mbox{where} \quad V_j=\left(\sum_{r=1}^n V_{rj}\right) \mbox{ mod }2.
\ena
In other words, in all stages other than $m$ and $m+1$ the switch variables that produce $V$ are the ones
from the given collection ${\bf X}$. In stage $m$, the switch variable for all bulbs but bulb $B_m$, chosen uniformly over
all bulbs in that stage that were not toggled, are the ones given by ${\bf X}$. The switch variable for $B_m$ in stage $m$, however,
is set to $C_m$, which takes the values $0$ and $1$ equally likely.
Hence with probability $1/2$, one additional bulb in stage $m$ is toggled. Similarly, in stage $m+1$, the switch
variable of bulb $B_{m+1}$, uniformly selected from all the bulbs that were toggled in that stage, is no longer
toggled with probability $1/2$. Since ${\bf X}$ and ${\bf V}$ differ in at most two switches, we have
\bea \label{difVX.is.2}
|X-V| \le 2.
\ena

We now show some basic facts about the distribution of the switch variables ${\bf V}$, that as in
even case $EV=n/2$, and verify the variance formula (\ref{defmuV}). In the following, let
\bea \label{def:nab}
{\bf n}(a,b)=(1,\ldots,m-1,m+a,m+b,m+2,\ldots,n) \qmq{for $a,b \in \{0,1\}$.}
\ena

\begin{lemma}
\label{lem:Vrdist}
The collections of variables $\{V_{rj},j=1,\ldots,n\}$ are mutually independent for $r=1,\ldots,n$.
For $r \in \{m,m+1\}$,
\bea \label{Vismix}
{\cal L}(V_{r1},\ldots,V_{rn})=\frac{1}{2}{\cal L}(X_{m1},\ldots,X_{mn})+\frac{1}{2}{\cal L}(X_{m+1,1},\ldots,X_{m+1,n}),
\ena
and furthermore, with $V_j$ given by (\ref{defVVk}), we have $P(V_j=1)=1/2$ for all $j=1,\ldots,n$, $EV=n/2$ and
$\mbox{Var}(V)=\sigma_V^2$ as given in (\ref{defmuV}).
\end{lemma}

\noindent \proof  The first claim follows by the independence of $\{X_{rj},j=1,\ldots,n\}$ for $r=1,\ldots,n$, and that $\{V_{rj},j=1,\ldots,n\}=\{X_{rj},j=1,\ldots,n\}$ for $r \not \in \{m,m+1\}$, and otherwise
$\{V_{rj},j=1,\ldots,n\}$ is given by randomizing $\{X_{rj},j=1,\ldots,n\}$ independently of the remaining stages,
and of stage $2m+1-r$.

To prove (\ref{Vismix}), first consider the case $r=m$. Let $e_1,\ldots,e_n \in \{0,1\}$. Since $(V_{m1},\ldots,V_{mn})=(X_{m1},\ldots,X_{mn})$ when $C_m=0$,
and since $C_m$ is independent of $(X_{m1},\ldots,X_{mn})$, we have
\beas
\lefteqn{P(V_{m1}=e_1,\ldots,V_{mn}=e_n)}\\
&=&\frac{1}{2}P(V_{m1}=e_1,\ldots,V_{mn}=e_n|C_m=0)+\frac{1}{2}P(V_{m1}=e_1,\ldots,V_{mn}=e_n|C_m=1)\\
&=&\frac{1}{2}P(X_{m1}=e_1,\ldots,X_{mn}=e_n|C_m=0)+\frac{1}{2}P(V_{m1}=e_1,\ldots,V_{mn}=e_n|C_m=1)\\
&=&\frac{1}{2}P(X_{m1}=e_1,\ldots,X_{mn}=e_n)+\frac{1}{2}P(V_{m1}=e_1,\ldots,V_{mn}=e_n|C_m=1).
\enas
Hence (\ref{Vismix}) is equivalent to
\bea \label{VmXmp1}
P(V_{m1}=e_1,\ldots,V_{mn}=e_n|C_m=1)=P(X_{m+1,1}=e_1,\ldots,X_{m+1,n}=e_n).
\ena
Now assuming $\sum_j e_j=m+1$, as both sides of (\ref{VmXmp1}) are zero otherwise, using the independence of
$C_m$ from $B_m$ and ${\bf X}$ we have
\beas
\lefteqn{P(V_{m1}=e_1,\ldots,V_{mn}=e_n|C_m=1)}\\
&=&\sum_{j=1}^n P(V_{m1}=e_1,\ldots,V_{mn}=e_n|B_m=j, C_m=1)P(B_m=j|C_m=1)\\
&=&\sum_{j:e_j=1}  P(X_{ml}=e_l,l \not = j, X_{mj}=0|B_m=j, C_m=1)P(B_m=j)\\
&=&\sum_{j:e_j=1} P(X_{ml}=e_l,l \not = j, X_{mj}=0|B_m=j)P(B_m=j)\\
&=&\sum_{j:e_j=1} P(B_m=j|X_{ml}=e_l,l \not = j, X_{mj}=0)P(X_{ml}=e_l,l \not = j, X_{mj}=0).
\enas
Now, since $B_m$ is chosen uniformly from the $n-m$ variables taking the value 0 in stage $m$, recalling $\sum_j e_j = m+1$, we have
\beas
\lefteqn{P(V_{m1}=e_1,\ldots,V_{mn}=e_n|C_m=1)=\frac{1}{n-m}\sum_{j:e_j=1} P(X_{ml}=e_l,l \not =j, X_{mj}=0)}\\
&=&\frac{m+1}{n-m} {n \choose m}^{-1} = {n \choose m+1}^{-1}=P(X_{m+1,1}=e_1,\ldots,X_{m+1,n}=e_n),
\enas
proving (\ref{VmXmp1}). The argument for $r=m+1$ is essentially the same.

Recalling the notation in (\ref{def:nab}), by (\ref{Vismix}) the distribution of ${\bf V}$ is the equal mixture of ${\cal L}({\bf X}_{\bf s})$ over the four cases ${\bf s}={\bf n}(a,b)$ for  $\{a,b\} \in \{0,1\}$. Hence, from (\ref{genmeanvar}),
\beas
EV=\frac{n}{8}\sum_{r,t \in \{m,m+1\}} \left( 1- \lambda_{n,1,{\bf n}(a,b)} \right) = \frac{n}{2},
\enas
as
\beas
\frac{1}{4}\sum_{a,b \in \{0,1\}} \lambda_{n,1,{\bf n}(a,b)}
&=&\frac{1}{4}\left( \prod_{s \not \in \{m,m+1\}}\lambda_{n,1,s} \right) \left(\lambda_{n,1,m}^2+ 2\lambda_{n,1,m}\lambda_{n,1,m+1}+ \lambda_{n,1,m+1}^2 \right)\\
&=&\left( \prod_{s \not \in \{m,m+1\}}\lambda_{n,1,s} \right) \left((\lambda_{n,1,m}+ \lambda_{n,1,m+1})/2 \right)^2=0,
\enas
since
\beas
\lambda_{n,1,m}+ \lambda_{n,1,m+1} = \left( 1-\frac{2m}{n} \right) +\left( 1-\frac{2(m+1)}{n} \right)=0.
\enas
Clearly $EV=n/2$ implies $P(V_j=1) = 1/2$

Next, by the conditional variance formula and (\ref{genmeanvar}) we obtain
\beas
\mbox{Var}(V) = \frac{1}{4}\sum_{a,b \in \{0,1\}}\left( \mbox{Var}(X_{{\bf n}(a,b)})+ E\left(X_{{\bf n}(a,b)}-\frac{n}{2}\right)^2 \right)
=\frac{1}{4}\sum_{a,b \in \{0,1\}}\left(\frac{n}{4}(1-\lambda_{n,2,{\bf n}(a,b)})+
\frac{n^2}{4}\lambda_{n,2,{\bf n}(a,b)}\right),
\enas
and now applying $(1/4)\sum_{a,b \in \{0,1\}} \lambda_{n,1,{\bf n}(a,b)} = \overline{\lambda}_{n, 1, {\bf n}}$, and a similar identity for $\lambda_{n,2,{\bf n}(a,b)}$,
yields the variance $\sigma_V^2$ in (\ref{defmuV}).
\bbox \\

We now present a size bias coupling for $V$. As in the even case, the variable $V^s$ is obtained by first constructing,
for each $i=1,\ldots,n$, switch variables ${\bf V}^i$ that satisfy
\bea \label{bfVi.correct.cdist}
{\cal L}({\bf V}^i)={\cal L}({\bf V}|V_i=1).
\ena
For a given $i=1,\ldots,n$, to construct ${\bf V}^i$, one first determines if $V_i=1$. If so, set ${\bf V}^i={\bf V}$. Otherwise, let $M$ be a variable that chooses from the stages $m$ and $m+1$ uniformly and independently of ${\bf V}$.
As in this case $V_i=0$, one may achieve $V_i^i=1$ by changing the switch variable $V_{Mi}$ to $1-V_{Mi}$. The coupling
accomplishes this change in one of two possible ways.

To introduce the first way, called a flip, we say a configuration
${\bf e}$ of binary switch variable values is feasible if $P({\bf V}={\bf e}) \not = 0$, that is, when
\beas 
\sum_{j=1}^n e_{rj} = r \qmq{for $r \not \in \{m,m+1\}$\,\,\,\,\,\,\, and}
\,\, \sum_{j=1}^n e_{rj} \in \{m,m+1\} \qmq{for $r \in \{m,m+1\}$.}
\enas
If flipping $V_{Mi}$ to $1-V_{Mi}$ results in feasible configuration, then the flip is made with probability $1/(m+1)$. In other words, given ${\bf e}$, $r \in \{m,m+1\}$ and $i=1,\ldots,n$ let ${\bf e}^{r,i}$ be the configuration with entries
\beas 
e_{sl}^{r,i}=\left\{
\begin{array}{cc}
e_{sl} & \mbox{$s \not =r$ or $l \not = i$}\\
1-e_{ri} & \mbox{$s=r$ and $l=i$}.
\end{array}
\right.
\enas
Defining ${\bf V}^{M,i}$ in like manner, the distribution of $F^i$, the indicator that $V_{Mi}$ is flipped, is given by
\bea \label{distFi}
P(F^i=1|{\bf V},M)=\frac{1}{m+1}{\bf 1}(\mbox{${\bf V}^{M,i}$ is feasible}).
\ena

In the flip is unsuccessful, that is, if $F^i=0$, we perform an `interchange' in stage $M$, much like the coupling in the even case. For a configuration ${\bf e}, r \in \{m,m+1\}$ and $i,j
\in \{1,\ldots,n\}$, let
\beas
e_{sl}^{r,i,j}=\left\{
\begin{array}{cc}
e_{sl} & \mbox{$s \not =r$ or $l \not = i$}\\
e_{rl} & s=r,l \not \in \{i,j\}\\
e_{ri} & s=r, l=j\\
e_{rj} & s=r, l=i,
\end{array}
\right.
\enas
that is, ${\bf e}^{r,i,j}$ is the configuration ${\bf e}$ with the variables in the $r,i$ and $r,j$ positions
interchanged.  Now let $J^i$ be a random index with
distribution given by
\bea \label{distJi}
{\cal L}(J^i|{\bf V},M) = {\cal U}\{j:V_{Mj} \not = V_{Mi}\}.
\ena
Defining ${\bf V}^{M,i,j}$ in like manner, when $F^i=0$, we interchange $V_{Mi}$ with $V_{M,J^i}$.

Hence, overall the configuration ${\bf V}^i$ is specified by
\bea \label{defVi.flip.interchange}
{\bf V}^i = \left\{
\begin{array}{ll}
{\bf V}       & V_i=1\\
{\bf V}^{M,i} & V_i=0, F^i=1\\
{\bf V}^{M,i,J^i} & V_i=0, F^i=0.
\end{array}
\right.
\ena
The following theorem shows that ${\bf V}^i$ satisfies (\ref{bfVi.correct.cdist}).

Before stating the theorem, we claim that if ${\bf e}$ is feasible,
then for $r \in\{m,m+1\}$, the configuration ${\bf e}^{r,i}$ is
feasible if and only if $\sum_{j \not = i} e_{rj}=m$. Both ${\bf e}$
and ${\bf e}^{r,i}$ are feasible if and only if $\sum_{j \not = i} e_{rj}
+ e_{ri} \in \{m,m+1\}$ and $\sum_{j \not = i} e_{rj} + 1- e_{ri} \in
\{m,m+1\}$. The claim now follows noting that
$\{m-e_{ri},m+1-e_{ri}\} \cap \{m-1+e_{ri}, m+e_{ri}\} = \{m-1,m\} \cap \{m,m+1\}=\{m\}$.
Summarizing, for a feasible configuration ${\bf e}$,
\bea \label{bferi.feasible.b}
{\bf e}^{r,i} \qmq{is feasible if and only if} \sum_{j \not = i} e_{rj}=m \quad \mbox{that is, if and only if $e_{ri}=\sum_{j=1}^n e_{rj}-m$.}
\ena

In the following we denote $F^i$ and $J^i$ by $F$ and $J$ respectively for simplicity.

\begin{theorem}
\label{thm:oddcop}
Let ${\bf V}$ be constructed from ${\bf X}$ as in (\ref{def:VX}), let $M$ be a random variable uniformly distributed over $\{m,m+1\}$ and independent ${\bf V}$, and for $i \in \{1,\ldots,n\}$, given ${\bf V}$ and $M$ let $F$ and $J$ have distributions as
specified in (\ref{distFi}) and (\ref{distJi}), respectively. Then ${\bf V}^i$ given by (\ref{defVi.flip.interchange}) satisfies (\ref{bfVi.correct.cdist}).

Further, letting $V^i=\sum_{j=1}^n V_j^i$ where
$$
V_j^i = \left(\sum_{r=1}^n V_{rj}^i\right) \mbox{ mod }2,
$$
when $I$ is uniformly
chosen from $\{1,\ldots,n\}$ independent of all other variables, the mixture $V^I=V^s$ has the $V$-size
bias distribution and satisfies
\bea \label{VVsbdmn}
V^s-V={\bf 1}_{\{V_I=0,F=1\}}+2 {\bf 1}_{\{V_I=0,
X_J=0,F=0\}} \qmq{and} V \le V^s \le V+2.
\ena

\end{theorem}

\noindent \proof  As in the even case, since ${\bf V}^i={\bf V}$ when $V_i=1$ and $P(V_i=1)=1/2$, the proof that $V^I$ has the $V$-size bias distribution follows by verifying, for all $i=1,\ldots,n$, that
\bea \label{condVi.is.V}
P({\bf V}^i={\bf e}|V_i=0)=P({\bf V}={\bf e}|V_i=1) \qmq{or equivalently} P({\bf V}^i={\bf e},V_i=0)=P({\bf V}={\bf e},V_i=1)
\ena
holds for all feasible configurations ${\bf e}$ satisfying $\sum_{r=1}^n e_{ri} = 1 \,\mbox{mod}\,2$.

Recalling that $M,F$ and $J$ are the stage selected, the indicator that the bit in question in the ${\bf V}$ configuration is flipped, and the index with which an interchange is to be performed, respectively, we find
\bea
\lefteqn{P({\bf V}^i={\bf e},V_i=0)} \nn \\
&=&  P({\bf V}^i={\bf e},V_i=0,F=1)+P({\bf V}^i={\bf e},V_i=0,F=0)\nn \\
&=& \sum_{r=m}^{m+1} \left( P({\bf V}^i={\bf e},V_i=0,M=r,F=1) + \sum_{j=1}^n  P({\bf V}^i={\bf e},V_i=0, J=j, M=r,F=0) \right) \nn \\
&=& \sum_{r=m}^{m+1} \left( P({\bf V}={\bf e}^{r,i},V_i=0,M=r,F=1) + \sum_{j=1}^n P({\bf V}={\bf e}^{r,i,j}, V_i=0, J=j, M=r,F=0) \right) \nn \\
&=& \sum_{r=m}^{m+1} \left( P({\bf V}={\bf e}^{r,i},M=r,F=1) + \sum_{j=1}^n P({\bf V}={\bf e}^{r,i,j},J=j, M=r,F=0) \right) \label{dec.by.F}
\ena
where in the final equality we have used that $\sum_{r=1}^n e_{ri} = 1 \,\mbox{mod}\,2$.

Let $\sum_{j=1}^n e_{rj}=m+b$, where necessarily $b \in \{0,1\}$.  For the first term in the $r^{th}$ summand of (\ref{dec.by.F}), by (\ref{bferi.feasible.b}), (\ref{distFi}), the independence of $M$ from ${\bf V}$ and that
$P({\bf V}={\bf e})$ has the same value for any feasible configuration, we have
\bea
\lefteqn{P({\bf V}={\bf e}^{r,i},M=r,F=1)}\nn \\
&=& P(F=1|{\bf V}={\bf e}^{r,i},M=r)P({\bf V}={\bf e}^{r,i},M=r)=
\frac{1}{2}P({\bf V}={\bf e})\frac{{\bf 1}(e_{ri}=b)}{m+1}. \label{first.r.summand}
\ena

Next, for the second term in the $r^{th}$ summand of (\ref{dec.by.F}), we have
\bea
\lefteqn{\sum_{j=1}^n P({\bf V}={\bf e}^{r,i,j},J=j, M=r,F=0)}\nn \\
&=&\sum_{j=1}^n P(J=j|F=0,{\bf V}={\bf e}^{r,i,j},M=r)P(F=0|{\bf V}={\bf e}^{r,i,j},M=r) P({\bf V}={\bf e}^{r,i,j},M=r)\nn \\
&=&\frac{1}{2} P({\bf V}={\bf e})  \sum_{j=1}^n P(J=j|F=0,{\bf V}={\bf e}^{r,i,j},M=r)P(F=0|{\bf V}={\bf e}^{r,i,j},M=r),
\label{2ndrthdec.by.F}
\ena
the final equality holding by again using that $M$ is uniformly chosen over $\{m,m+1\}$, independent of all other variables.

To handle the summand terms in (\ref{2ndrthdec.by.F}), note that the index $J$ chooses uniformly over the $m+b(1-e_{ri}^{r,i,j})+(1-b)e_{ri}^{r,i,j}$ indices whose stage $r$ switch value is opposite to $e_{ri}^{r,i,j}$, that is,
\bea \label{Pjmm}
P(J=j|F=0,{\bf V}={\bf e}^{r,i,j},M=r)=\frac{{\bf 1}(e_{ri}^{r,i,j} \not = e_{rj}^{r,i,j})}{m+b(1-e_{ri}^{r,i,j})+(1-b)e_{ri}^{r,i,j}}=\frac{{\bf 1}(e_{rj}\not = e_{ri})}{m+b e_{ri}+(1-b)(1-e_{ri})}.
\ena
Again using $e_{ri}^{r,i,j}=e_{rj}$, by (\ref{bferi.feasible.b}) and (\ref{distFi}),
\bea \label{PFmm}
P(F=0|{\bf V}={\bf e}^{r,i,j},M=r)
={\bf 1}(e_{rj}=1-b) + \frac{m}{m+1}{\bf 1}(e_{rj}=b).
\ena
Taking the product of (\ref{Pjmm}) and (\ref{PFmm}) yields
\bea \label{49*48=2352}
\frac{{\bf 1}(e_{rj}\not = e_{ri})}{m+b e_{ri}+(1-b)(1-e_{ri})} \left({\bf 1}(e_{rj}=1-b) + \frac{m}{m+1}{\bf 1}(e_{rj}=b)\right)=
\frac{1}{m+1}{\bf 1}(e_{ri} \not = e_{rj}).
\ena
Hence (\ref{2ndrthdec.by.F}) may be written as
\beas
\frac{P({\bf V}={\bf e})}{2(m+1)} \sum_{j=1}^n {\bf 1}(e_{ri} \not = e_{rj})=
\frac{1}{2}P({\bf V}={\bf e})\left( \frac{m}{m+1}{\bf 1}(e_{ri}=b) + {\bf 1}(e_{ri}=1-b) \right).
\enas
Adding this factor to (\ref{first.r.summand}) we find that the $r^{th}$ summand of (\ref{dec.by.F}) equals
\beas
\frac{1}{2}P({\bf V}={\bf e})\frac{{\bf 1}(e_{ri}=b)}{m+1}+ \frac{1}{2}P({\bf V}={\bf e})\left( \frac{m}{m+1}{\bf 1}(e_{ri}=b) + {\bf 1}(e_{ri}=1-b) \right)=\frac{1}{2}P({\bf V}={\bf e}).
\enas
Hence, adding up the two summands $r=m$ and $r=m+1$ in (\ref{dec.by.F}) yields (\ref{condVi.is.V}), as
\beas
P({\bf V}^i={\bf e},V_i=0) = P({\bf V}={\bf e}) = P({\bf V}={\bf e},V_i=1),
\enas
where we have applied $\sum_{r=1}^n e_{ri} = 1 \,\mbox{mod}\,2$ for the final equality.

Lastly, to prove (\ref{VVsbdmn}), first note that by construction $V^s=V$ when $V_I=1$. When $V_I=0$ and $F=0$ then $V^s$ equals either
$V$ or $V+2$, depending on whether $X_J$ takes the value $1$ or $0$, respectively, as in the even case. When $V_I=0$ and $F=1$ the status of bulb $I$ is changed from off to on and the status of all other bulbs remain unchanged, hence $V^s-V=1$. \bbox

With the coupling of $V$ to $V^s$ now in hand, we prove a bound to the normal for $V$. We recall from Lemma \ref{lem:Vrdist}
that $EV=n/2$ and $\mbox{Var}(V)=\sigma_V^2$, given in (\ref{defmuV}).
\begin{theorem}
\label{thm:normV}
If $n=2m+1$, an odd number, and $V$ is given by (\ref{defVVk}), then
\beas
\sup_{z \in
\mathbb{R}} \left|P\left( \frac{V-n/2}{\sigma_V} \le z \right)-P(Z \le
z)\right| \le \frac{n}{2\sigma_V^2}{\overline \Delta}_1 + 1.64 \frac{n}{\sigma_V^3}+
\frac{2}{\sigma_V} \qmq{for all $n \ge 7$,}
\enas
where $\sigma_V^2$ and $\overline{\Delta}_1$
are given in (\ref{defmuV}) and (\ref{upDelta1}) respectively.
\end{theorem}

\noindent \proof  We apply Theorem \ref{L-BD:SDKsConc} to the coupling construction given in Theorem \ref{thm:oddcop}. As $\delta=2/\sigma_V$ from (\ref{VVsbdmn}), and $EV=n/2$ from Lemma \ref{lem:Vrdist}, it remains only to show that $\sqrt{\mbox{Var}\left( E\left( V^s-V \vert V \right)\right)} \le {\overline \Delta}_1$ for all $n \ge 7$. As in the even case, we obtain an upper bound by
conditioning on the $\sigma$-algebra ${\cal F}$ generated by ${\bf V}$, with respect to which $V$ is measurable.
Again for notational simplicity we drop the superscripts on $F$ and $J$ unless clarity demands them.

Taking conditional expectation in (\ref{VVsbdmn}) yields
\bea \label{def:zetan.xin}
E(V^s-V|{\cal F}) = \zeta_n + 2\xi_n  \qmq{where} \zeta_n= E\left( {\bf 1}_{\{V_I=0,F=1\}}|{\cal F}\right) \qmq{and}
\xi_n = E\left( {\bf 1}_{\{V_I=0, V_J=0,F=0\}}|{\cal F}\right).
\ena
Letting
\beas
A= \mbox{Var}\left( \zeta_n \right) \qmq{and}
B= \mbox{Var} \left( \xi_n \right),
\enas
since $|\mbox{Cov}(\zeta_n,\xi_n)| \le \sqrt{AB}$ by the Cauchy-Schwarz inequality yields, we have
\bea \label{improved.CS.bound}
\mbox{Var}(E(V^s-V|{\cal F})) \le A+4\sqrt{AB}+4B.
\ena

First computing a bound on $A$, expanding the indicator defining $\zeta_n$ in (\ref{def:zetan.xin}) over the possible values of $M$ and $I$ yields
\beas
{\bf 1}_{\{V_I=0,F=1\}}
= \sum_{r=m}^{m+1} \sum_{i=1}^n    {\bf 1}_{\{M=r,I=i,V_i=0,F=1 \}}.
\enas
Recall that conditional on ${\bf V}$, configurations
for which ${\bf V}^{r,i}$ are feasible, that is, by (\ref{bferi.feasible.b}), those for which $\sum_{j \not =i}V_{rj}=m$, are flipped with probability $1/(m+1)$. Now, since
$I$ and $M$ are chosen uniformly and independently over indices $\{1,\ldots,n\}$ and $\{m,m+1\}$ respectively, taking conditional expectation
with respect to ${\cal F}$ we obtain
\bea \label{exp:term 1}
\zeta_n
= \frac{1}{2n(m+1)}\sum_{r=m}^{m+1}   \sum_{i=1}^n
{\bf 1}\left\{V_i=0,
\sum_{j \not =i}V_{rj}=m
\right\}.
\ena

To determine $E\zeta_n$ we now demonstrate that the events $\{V_i = 0\}$ and $\{\sum_{j \not =i}V_{rj}=m\}$ are independent for $r \in \{m,m+1\}$. Note that for each fixed $i$,
\begin{eqnarray*}
P\left(\sum_{j \not =i}V_{rj}=m\right) = P\left(V_i = 0, \sum_{j \not =i}V_{rj}=m\right)+P\left(V_i = 1, \sum_{j \not =i}V_{rj}=m\right),
\end{eqnarray*}
so it suffices to prove that
\bea \label{PV1.PV0}
P\left(V_i = 0, \sum_{j \not =i}V_{rj}=m\right) = P\left(V_i = 1, \sum_{j \not =i}V_{rj}=m\right)
\ena
since $P(V_i = 0) = 1/2$, by Lemma \ref{lem:Vrdist}. However, as the map that sends ${\bf e}$ to ${\bf e}^{r,i}$ is a bijection between the two events in (\ref{PV1.PV0}), since all feasible configurations have equal probability the equality in (\ref{PV1.PV0}) holds, and $\{V_i = 0\}$ and $\{\sum_{j \not =i}V_{rj}=m\}$ are therefore independent.

Now, using that
\bea
P\left(\sum_{j\ne i}V_{rj} = m\right) &=& P\left(\sum_{j\ne i}V_{rj} = m, C_r = 0\right) + P\left(\sum_{j\ne i}V_{rj} = m, C_r = 1\right) \nn \\
&=& \frac{1}{2} P\left(\sum_{j\ne i}V_{rj} = m \Bvert C_r = 0\right) + \frac{1}{2} P\left(\sum_{j\ne i}V_{rj} = m \Bvert C_r = 1\right) \nn \\
&=& \frac{1}{2} P(V_{ri}=0|C_r=0)+\frac{1}{2}P(V_{ri}=1|C_r=1) \nn \\
&=& \frac{1}{2} \left( \frac{m+1}{n} +\frac{m+1}{n} \right) = \frac{m+1}{n}, \label{Vmj=m}
\ena
and that $P(V_i=0)=1/2$, taking expectation and using independence in (\ref{exp:term 1}) we have
\bea \label{EUn}
E\zeta_n = \frac{1}{2n(m+1)}\sum_{r=m}^{m+1} \sum_{i=1}^n \frac{m+1}{2n} = \frac{1}{2n}.
\ena

Turning to the second moment of $\zeta_n$, squaring the right hand side of (\ref{exp:term 1}) we obtain
\[
\zeta_n^2 = \frac{1}{4n^2(m+1)^2}
\sum_{r,t \in \{m,m+1\}}\sum_{i=1}^n \sum_{l = 1}^n {\bf 1}\left\{V_i = 0, V_l = 0, \sum_{j \not =i}V_{rj}=m, \sum_{j \not =l}V_{tj}=m\right\}.
\]

We decompose $4n^2(m+1)^2\zeta_n^2$ into the following four terms:
\begin{enumerate}
\item For $r=t,i=l$ we obtain $\zeta_{n,1}^2 = \sum_{r=m}^{m+1} \sum_{i=1}^n {\bf 1}\{V_i = 0, \sum_{j \not =i}V_{rj}=m\}$

\item For $r=t,i\not =l$ we obtain  $\zeta_{n,2}^2 = \sum_{r=m}^{m+1}\sum_{i \not =l} {\bf 1}\{V_i = 0, V_l = 0, \sum_{j \not =i}V_{rj}=m, \sum_{j \not =l}V_{rj}=m\}$

\item For $r\not =t,i=l$ we obtain $\zeta_{n,3}^2 = \sum_{r \not = t} \sum_{i=1}^n
{\bf 1}\{V_i = 0, \sum_{j \not =i}V_{rj}=m, \sum_{j \not =i}V_{tj}=m\}$

\item For $r\not =t,i\not =l$ we obtain $\zeta_{n,4}^2 = \sum_{r \not = t} \sum_{i \not = l}
{\bf 1}\{V_i = 0, V_l=0, \sum_{j \not =i}V_{rj}=m, \sum_{j \not =l}V_{tj}=m\}$

\end{enumerate}

For the expectation of $\zeta_{n,1}^2$, using independence as in the
calculation of $E\zeta_n$ and (\ref{Vmj=m}), we obtain
\bea \label{Un12}
E\zeta_{n,1}^2 = \sum_{r=m}^{m+1} \sum_{i=1}^n P(V_i = 0, \sum_{j \not =i}V_{rj}=m) = m+1.
\ena

Turning to $\zeta_{n,2}^2$ we have
\beas
&&\{V_i = 0, V_l = 0, \sum_{j \not = i} V_{rj}=m, \sum_{j \not = l} V_{rj}=m\}\\
&&=\bigcup_{b \in \{0,1\}} \{V_i = 0, V_l = 0, \sum_{j \not \in
\{i,l\}} V_{rj}=m-b,\, V_{ri}=b, V_{rl}=b\},
\enas
and hence may write
\bea \label{Un2sq}
\zeta_{n,2}^2=\sum_{b \in \{0,1\}} \zeta_{n,2,b}^2 \qmq{where}
\zeta_{n,2,b}^2=\sum_{r=m}^{m+1}\sum_{i \not = l} {\bf 1}\left\{V_i = 0,
V_l = 0, \sum_{j \not \in \{i,l\}} V_{rj}=m-b,\, V_{ri}=b,
V_{rl}=b\right\}.
\ena
Hence, to compute $E\zeta_{n,2,b}^2$ we sum
\beas
&&P\left( V_i = 0, V_l = 0, \sum_{j \not \in \{i,l\}} V_{rj}=m-b, \,
V_{ri}=b, V_{rl}=b \right)\\
&&=P\left( V_i = 0, V_l = 0, \sum_{j \not \in \{i,l\}} V_{mj}=m-b, \,
V_{mi}=b, V_{ml}=b \right)\\
&&= P\left( V_i = 0, V_l = 0, \sum_{j=1}^n V_{mj}=m+b, \, V_{mi}=b,V_{ml}=b \right)\\
&&= P\left( V_i = 0, V_l = 0, V_{mi}=b,V_{ml}=b, C_m=b \right)\\
&&= \sum_{a \in \{0,1\}} P\left( V_i = 0, V_l = 0, V_{mi}=b,V_{ml}=b, C_m=b, C_{m+1}=a \right)\\
&&= \frac{1}{4}\sum_{a \in \{0,1\}} P\left( V_1 = 0, V_2 = 0, V_{m1}=b,V_{m2}=b | C_m=b, C_{m+1}=a
\right),
\enas
where in the second line we used the exchangeability of stages $m$ and $m+1$ in ${\bf V}$, and for
the final equality that $i \not = l$, and exchangeability again.

Recalling definition (\ref{def:nab}), for $a,b, \in \{0,1\}$ easily we have that
\beas
{\cal L}({\bf V}|C_m=b,C_{m+1}=a)={\cal L}({\bf X}_{{\bf n}(b,a)}),
\enas
and so, applying definition (\ref{def:galphabeta}), we obtain
\beas
P\left( V_1 = 0, V_2 = 0, V_{m1}=b,V_{m2}=b | C_m=b, C_{m+1}=a
\right) = g_{2-2b,2b,{\bf n}(b,a)}^{(m)}.
\enas
Hence, recalling (\ref{Un2sq}),
\bea \label{Un22}
E\zeta_{n,2}^2 = \frac{(n)_2}{2} \sum_{a,b \in \{0,1\}} g_{2-2b,2b,{\bf n}(b,a)}^{(m)}.
\ena

Now turning to $E\zeta_{n, 3}^2$, by again considering the bijection that maps ${\bf e}$ to ${\bf e}^{r,i}$ we may conclude that
\beas
P\left( V_i = 0, \sum_{j \not =i}V_{rj}=m, \sum_{j \not =i}V_{tj}=m \right)=P\left( V_i = 1, \sum_{j \not =i}V_{rj}=m, \sum_{j \not =i}V_{tj}=m \right)
\enas
and hence, by the independence of the switch variables in stages $r \not = t$, as provided by Lemma \ref{lem:Vrdist}, we obtain
\beas
\frac{(m+1)^2}{n^2} = P\left(\sum_{j \not =i}V_{rj}=m, \sum_{j \not =i}V_{tj}=m \right) = 2P\left( V_i = 0, \sum_{j \not =i}V_{rj}=m, \sum_{j \not =i}V_{tj}=m \right),
\enas
and therefore
\bea \label{Un32}
E\zeta_{n,3}^2 = \frac{(m+1)^2}{n}.
\ena

Arguing similarly to compute $\zeta_{n,4}^2$, the identity
\beas
\frac{(m+1)^2}{n^2}= P\left(\sum_{j \not =i}V_{rj}=m, \sum_{j \not =l}V_{tj}=m \right)  = 4 P\left( V_i = 0, V_l=0, \sum_{j \not =i}V_{rj}=m, \sum_{j \not =l}V_{tj}=m \right)
\enas
yields
\bea \label{Un42}
E\zeta_{n,4}^2=\frac{(n)_2 (m+1)^2}{2n^2}.
\ena

Combining (\ref{Un12}), (\ref{Un22}), (\ref{Un32}), (\ref{Un42}) and (\ref{EUn}) yields
\beas 
\lefteqn{A = E(\zeta_n^2) - E(\zeta_n)^2}\\
&=& \frac{1}{4n^2(m+1)^2}\left( (m+1) + \frac{(n)_2}{2} \sum_{a,b \in \{0,1\}} g_{2-2b,2b,{\bf n}(b,a)}^{(m)} + \frac{(m+1)^2}{n} + \frac{(n)_2(m+1)^2}{2n^2}\right) - \frac{1}{4n^2}.\nn
\enas

To obtain a bound on $A$ first note that for any $b$ we have ${\bf n}(b,0)_m={\bf n}_{m+1}$ and ${\bf n}(b,1)_m={\bf n}_m$, and that the $m^{th}$ component of ${\bf n}(b,a)$ is $m+b$. Hence when $b=0$ Corollary \ref{cor:spectral.4} yields

\beas
g_{2,0,{\bf n}(0,a)}^{(m)} = \left\{
\begin{array}{ll}
\frac{1}{4}(1+2\lambda_{n,1,{\bf n}_{m+1}}+ \lambda_{n,2,{\bf n}_{m+1}})\frac{(m+1)_2}{(n)_2} & a=0 \\
\frac{1}{4}(1+2\lambda_{n,1,{\bf n}_m}+ \lambda_{n,2,{\bf n}_m})\frac{(m+1)_2}{(n)_2} & a=1,
\end{array}
\right.
\enas
and when $b=1$ that
\beas
g_{0,2,{\bf n}(1,a)}^{(m)} = \left\{
\begin{array}{ll}
\frac{1}{4}(1-2\lambda_{n,1,{\bf n}_{m+1}}+ \lambda_{n,2,{\bf n}_{m+1}})\frac{(m+1)_2}{(n)_2} & a=0 \\
\frac{1}{4}(1-2\lambda_{n,1,{\bf n}_m}+ \lambda_{n,2,{\bf n}_m})\frac{(m+1)_2}{(n)_2} & a=1.
\end{array}
\right.
\enas
Summing these four terms, and using that $(\lambda_{n,2,{\bf n}_m}+\lambda_{n,2,{\bf n}_{m+1}})/2=\overline{\lambda}_{n,2,{\bf n}_m}$, as defined in (\ref{def:overlinelambda}), yields
\beas
\sum_{a, b \in \{0, 1\}} g_{2-2b, 2b,{\bf n}(b,a)}^{(m)}
&=& \frac{1}{4}\left( 4 + 2\lambda_{n,2,{\bf n}_{m}} +2 \lambda_{n,2,{\bf n}_{m+1}}\right)\frac{(m+1)_2}{(n)_2} = (1 + \overline{\lambda}_{n, 2, {\bf n}_m})\frac{(m+1)_2}{(n)_2}.
\enas
Hence,
\beas
\lefteqn{A}\\
&=& \frac{1}{4n^2(m+1)^2}\left( (m+1) + \frac{(n)_2}{2} \left( \frac{(m+1)_2}{(n)_2} + \overline{\lambda}_{n, 2, {\bf n}_m}\frac{(m+1)_2}{(n)_2} \right) + \frac{(m+1)^2}{n} + \frac{(n)_2}{2}\frac{(m+1)^2}{n^2}\right) - \frac{1}{4n^2}\nn \\
&=& \frac{1}{4n^2} \left(\frac{1}{m+1}+\frac{m}{2(m+1)}+\frac{1}{n}+\frac{n-1}{2n} - 1 \right) + \frac{m}{8n^2(m+1)}\overline{\lambda}_{n, 2, {\bf n}_m} \nn \\
&=& \frac{3m+2}{8n^3(m+1)} + \frac{m}{8n^2(m+1)}\overline{\lambda}_{n, 2, {\bf n}_m}. 
\enas
Applying Lemma \ref{lam24bds} we obtain
\bea \label{final.A}
A \le \frac{3m+2}{8n^3(m+1)} + \frac{e^{-n}}{8n^2} \le \frac{3}{8n^3} + \frac{e^{-n}}{8n^2}.
\ena

Next, expanding the indicator in (\ref{def:zetan.xin}) to obtain $B=\mbox{Var}(\xi_n)$, and recalling that an interchange is performed when $F=0$, we have
\beas
\lefteqn{{\bf 1}_{\{V_I=0, V_J=0,F=0\}} }\\
&=& \sum_{r=m}^{m+1} \sum_{i,j=1}^n {\bf 1}_{\{V_i=0, V_j=0 \}} {\bf 1}_{\{I=i,J=j,M=r,F=0\}} \\
&=& \sum_{r=m}^{m+1} \sum_{i,j=1}^n \left( {\bf 1}_{\{V_i=0, V_j=0, V_{ri}=0\}}+{\bf 1}_{\{V_i=0, V_j=0, V_{ri}=1\}}\right) {\bf 1}_{\{I=i,J=j,M=r,F=0\}}\\
&=& \sum_{r=m}^{m+1} \sum_{i \not = j} \left(  {\bf 1}_{\{V_i=0, V_j=0, V_{ri}=0,V_{rj}=1\}} + {\bf 1}_{\{V_i=0, V_j=0, V_{ri}=1,V_{rj}=0\}}\right) {\bf 1}_{\{I=i,J=j,M=r,F=0\}} \\
&=&  2\sum_{r=m}^{m+1}\sum_{i \not = j} {\bf 1}_{\{V_i=0, V_j=0, V_{ri}=0,V_{rj}=1\}}{\bf 1}_{\{I=i,J=j,M=r,F=0\}},
\enas
where the last inequality holds due to the symmetry between $i$ and $j$ in the penultimate sum, as in the even case.

As ${\bf V}$ is ${\cal F}$-measurable, recalling the definition $\xi_n=E({\bf 1}_{\{V_I=0,V_J=0,F=0\}}|{\cal F})$ in (\ref{def:zetan.xin}), taking conditional expectation yields
\bea \label{defxin.odd}
\xi_n = 2\sum_{r=m}^{m+1} \sum_{i\ne j}{\bf 1}_{\{V_i=0, V_j=0, V_{ri}=0,V_{rj}=1\}} E({\bf 1}_{\{I=i,J=j,M=r,F=0\}}|{\cal F}).
\ena
To compute the conditional expectation in the sum, using the
independence of $I$ from $M$ and ${\bf V}$, and $M$ from ${\bf V}$,
we obtain
\beas
\lefteqn{E({\bf 1}_{\{I=i,J=j,M=r,F=0\}}|{\cal F})}\\
&=& P(J=j|M=r,I = i, F=0,{\bf V}) P(F=0|I =i, M=r, {\bf V})P(I = i|M = r, {\bf V})P(M = r| {\bf V}) \\
&=& \frac{1}{2n} P(J=j|M=r,I = i, F=0,{\bf V}) P(F=0|I =i, M=r, {\bf V}) \\
&=& \frac{1}{2n} P(J^i=j|M=r,F^i=0,{\bf V}) P(F^i=0|M=r, {\bf V}) \\
&=& \frac{1}{2n}  \frac{{\bf 1}\{V_{ri} \ne V_{rj}\}}{m+1},
\enas
the final equality since the product of (\ref{PFmm}) and (\ref{Pjmm}) is given by (\ref{49*48=2352}).
Hence, as the first indicator in the sum in (\ref{defxin.odd}) specifies that $V_{ri}=0,V_{rj}=1$, we obtain
\bea \label{exp:xi}
\xi_n = \frac{1}{n(m+1)} \sum_{r=m}^{m+1} \sum_{i\ne j}{\bf 1}_{\{V_i=0, V_j=0, V_{ri}=0,V_{rj}=1\}}.
\ena

Taking the expectation of $\xi_n$, by the exchangeability in ${\bf V}$ of stages $m$ and $m+1$, we have
\bea \label{Exin}
E(\xi_n) = \frac{2}{n(m+1)}\sum_{i \not = j}P(V_i=0, V_j=0, V_{mi}=0,V_{mj}=1).
\ena
Exhausting over the possible values of $C_m = b$ and
$C_{m+1} = a$ for $a, b \in \{0, 1\}$, we obtain
\beas
\lefteqn{P(V_i=0, V_j=0, V_{mi}=0,V_{mj}=1)}\\
&=& \frac{1}{4}\sum_{a,b \in \{0, 1\}} P(V_i=0, V_j=0, V_{mi}=0,V_{mj}=1|C_m = b, C_{m+1} = a)\\
&=& \frac{1}{4}\sum_{a,b \in \{0, 1\}} g_{1,1,{\bf n}(b, a)}^{(m)}.
\enas
Now, by (\ref{Exin}), using Corollary \ref{cor:spectral.4},
\bea \label{mean.xi.n}
E \xi_n = \frac{n-1}{2(m+1)}\sum_{a, b \in \{0, 1\}} g_{1, 1,
{\bf n}(b, a)}^{(m)}  = \frac{n-1}{2(m+1)}\left(1 -
\overline{\lambda}_{n,2,{\bf n}_m} \right)\frac{m(m+1)}{(n)_2} =
\frac{m}{2n}\left(1 - \overline{\lambda}_{n,2,{\bf n}_m} \right).
\ena

To calculate $E\xi_n^2$, squaring $\xi_n$ in (\ref{exp:xi}) we obtain a sum over
indices $r,t \in \{m,m+1\}$ and $i_1,i_2,j_1,j_2$ with $i_1 \ne j_1,
i_2 \ne j_2$, and therefore $|\{i_1,i_2,j_1,j_2\}| =p$ for $p\in \{2, 3,
4\}$. Hence we may write
\begin{eqnarray*}
& & \xi_n^2 = \sum_{p \in \{1,2,3,4\}} \xi_{n,p}^2 \quad \mbox{where} \quad \xi_{n,p}^2=\sum_{q \in \{1,2\}}\xi_{n, p, q}^2 \quad \mbox{and}\\
& & \xi_{n, p, q}^2 = \frac{1}{n^2(m+1)^2} \sum_{|\{r,t\}|=q}
\sum_{\stackrel{|\{i_i,i_2,j_1,j_2\}|=p}{i_1 \ne j_1, i_2 \ne j_2}}{\bf 1}_{\{V_{i_1}=0, V_{j_1}=0,
V_{r,i_1}=0,V_{r,j_1}=1\}} {\bf 1}_{\{V_{i_2}=0, V_{j_2}=0,
V_{t,i_2}=0,V_{t,j_2}=1\}}.
\end{eqnarray*}

First consider the case $p=4$, that is, $\{i_1,i_2,j_1,j_2\}|=4$.
We begin with $\xi_{n, 4, 1}^2$, that is, when $r=t$. Using again that $\{C_r = b, C_{n-r} = a\}$ partitions
the space for $a, b \in \{0,1\}$, we have
\begin{eqnarray*}
\lefteqn{P(V_{i_1}=0,
V_{i_2}=0,V_{j_1}=0,V_{j_2}=0,V_{r,i_1}=0, V_{r,i_2}=0, V_{r,j_1}=1,
V_{r,j_2}=1, C_r = b, C_{n-r} = a)} \\
&=& \frac{1}{4} P(V_{i_1}=0,
V_{i_2}=0,V_{j_1}=0,V_{j_2}=0,V_{r,i_1}=0, V_{r,i_2}=0, V_{r,j_1}=1,
V_{r,j_2}=1|C_r = b, C_{n-r} = a) \\
&=& \frac{1}{4} P(V_1=0,
V_2=0,V_3=0,V_4=0,V_{r,1}=0, V_{r,2}=0, V_{r,3}=1,
V_{r,4}=1|C_r = b, C_{n-r} = a) \\
&=& \frac{1}{4} g_{2, 2,{\bf n}(b, a)}^{(r)}.
\end{eqnarray*}
Summing over the four distinct indices $i_1,i_2,j_1,j_2$ and applying Corollary \ref{cor:spectral.4}
we obtain
\bea
E(\xi_{n, 4, 1}^2)
&=&\frac{(n)_4}{4n^2(m+1)^2} \sum_{a,b \in \{0,1\}} \sum_{r=m}^{m+1} g_{2, 2,
{\bf n}(b, a)}^{(r)} \nn \\
&=& \frac{2(n)_4}{16n^2(m+1)^2}\left(1 - 2 \overline{\lambda}_{n,2,{\bf n}_{m}} + \overline{\lambda}_{n,4,{\bf n}_{m}}\right)\frac{(m)_2(m+1)_2}{(n)_4} \nn \\
&=&\frac{m^2(m-1)}{8n^2(m+1)}\left(1 - 2 \overline{\lambda}_{n,2,{\bf n}_{m}} + \overline{\lambda}_{n,4,{\bf n}_{m}}\right). \label{xin4.1}
\ena

In order to proceed further we need to consider functions such as those in (\ref{def:galphabeta}), but which specify switch variables
in two stages, rather than in only one. In particular, for ${\bf s}=(s_1,\ldots,s_k)$, distinct stages $r,t \in \{1,\ldots,k\}$ and bulbs $i_1<\cdots<i_\alpha$ and $j_1<\cdots<j_\beta$ in $\{1,\ldots,n\}$, let
\bea \label{def:g.2.stage}
\lefteqn{g_{[(i_1,a_1),\ldots,(i_\alpha,a_\alpha);(j_1,b_1),\ldots,(j_\beta,b_\beta)],{\bf
s}}^{(r,t)}}\\
&&=P(X_{i_1}=0,\ldots,X_{i_\alpha}=0,X_{j_1}=0,\ldots,X_{j_\beta}=0,X_{r,i_1}=a_1,\ldots,X_{r,i_\alpha}=a_\alpha,
X_{t,j_1}=b_1,\ldots,X_{t,j_\beta}=b_\beta) \nn
\ena
when ${\bf X}$ has switch pattern ${\bf s}$.
We note that it is not required that $\{i_1,\ldots,i_\alpha\} \cap  \{j_1,\ldots,j_\beta\} \not = \emptyset$.
The functions in (\ref{def:g.2.stage}) give the probability that bulbs $\{i_1,\ldots,i_\alpha\} \cup \{j_1,\ldots,j_\beta\}$ are off, and that the switch values for bulbs $i_1,\ldots,i_\alpha$ and $j_1,\ldots,j_\beta$
assume the prescribed values in stages $r$ and $t$, respectively.

Now moving to $\xi_{n, 4, 2}^2$, considering first $r=m,t=m+1$ and arguing similarly as for the $\xi_{n,4,1}$,
\begin{eqnarray*}
& & P(V_{i_1}=0,
V_{i_2}=0,V_{j_1}=0,V_{j_2}=0,V_{m,i_1}=0, V_{m+1,i_2}=0, V_{m,j_1}=1,
V_{m+1,j_2}=1, C_m = b, C_{m+1} = a) \\
&&= \frac{1}{4}P(V_1=0, V_2=0, V_3=0,V_4=0, V_{m,1}=0,
V_{m,2}=1, V_{m+1,3}=0,
V_{m+1,4}=1|C_m = b, C_{m+1} = a) \\
&&= \frac{1}{4}g_{[(1,0), (2, 1); (3, 0), (4, 1)], {\bf
n}(b, a)}^{(m,m+1)},
\end{eqnarray*}
applying the definition (\ref{def:g.2.stage}). Summing over the four distinct indices, noting that the case $r=m+1,t=m$ contributes equally, using (\ref{g:xi-42}) from the Appendix we obtain
\bea \nn
E(\xi_{n, 4, 2}^2) &=& \frac{(n)_4}{2n^2(m+1)^2} \sum_{a,b \in \{0,1\}} g_{[(1,0), (2, 1); (3, 0), (4, 1)],
{\bf n}(b, a)}^{(m,m+1)}\\
&=& \frac{(n)_4}{2n^2(m+1)^2}\left( \frac{(m+1)_3}{(n)_4}\right)^2 \left( (2m-1)^2 + 2(2m-1)\lambda_{n,2,{\bf n}_{m,m+1}} + \lambda_{n,4,{\bf n}_{m,m+1}} \right). \label{xi-n-4:exp}
\ena
Summing $E(\xi_{n, 4, 1}^2)$ and $E(\xi_{n, 4, 2}^2)$ from (\ref{xin4.1}) and (\ref{xi-n-4:exp}) respectively, and simplifying the latter using $n=2m+1$ yields
\bea \label{Exin42.sq}
E(\xi_{n, 4}^2)
&=& \frac{m^2(m-1)}{8n^2(m+1)}\left(1 - 2 \overline{\lambda}_{n,2,{\bf n}_{m}} + \overline{\lambda}_{n,4,{\bf n}_{m}}\right) \\
&&+ \frac{(m+1)_3^2}{2n^2(m+1)^2(n)_4} \left( (2m-1)^2  + 2\left(2m-1  \right) \lambda_{n,2,{\bf n}_{m,m+1}}  + \lambda_{n,4,{\bf n}_{m,m+1}} \right)\nn \\
&=&\frac{m^2(m-1)}{8n^2(m+1)}\left(1 - 2
\overline{\lambda}_{n,2,{\bf n}_{m}} + \overline{\lambda}_{n,4,{\bf
n}_{m}}\right) + \frac{(m)_2(2m-1)}{8n^3} +
\frac{(m)_2}{4n^3}\lambda_{n,2,{\bf n}_{m,m+1}}
+\frac{(m)_2^2}{2n^2(n)_4}\lambda_{n,4,{\bf n}_{m,m+1}}.\nn
\ena

Moving to $p=3$ and the evaluation of $E(\xi_{n,3,q}^2)$, under the constraints $i_1 \ne j_1, i_2 \ne j_2$, sums over indices that satisfy $|\{i_1,i_2,j_1,j_2\}|=3$
can be partitioned into the following four cases:
\beas
1.\, i_1=i_2, j_1 \not = j_2, \quad 2.\,i_1 \not = i_2, j_1=j_2, \quad 3.\, i_1=j_2, i_2 \not = j_1, \quad 4.\,i_2=j_1, i_1 \not = j_2,
\enas
Subscripting by these subcases, here and similarly in what follows, we write
\beas
\xi_{n,3,q}^2= \xi_{n,3,q,1}^2 + \xi_{n,3,q,2}^2+\xi_{n,3,q,3}^2 + \xi_{n,3,q,4}^2.
\enas

Starting with $q=1$, for $E\xi_{n,3,1,1}^2$, we have
\begin{eqnarray*}
\lefteqn{P(V_i=0 ,V_{j_1}=0,V_{j_2}=0,V_{r,i}=0, V_{r,j_1}=1, V_{r,j_2}=1, C_r = b, C_{n-r} = a)}\\
&=& \frac{1}{4}P(V_1=0 ,V_2=0,V_3=0,V_{r,1}=0, V_{r,2}=1, V_{r,3}=1|C_r = b, C_{n-r} = a) \\
&=& \frac{1}{4} g_{1, 2, {\bf n}(b, a)}^{(r)},
\end{eqnarray*}
while for $E\xi_{n,3,1,2}^2$,
\begin{eqnarray*}
\lefteqn{P(V_{i_1}=0 ,V_{i_2}=0,V_{j}=0,V_{r,i_1}=0, V_{r,i_2}=0, V_{r,j}=1, C_r = b, C_{n-r} = a)}\\
&=& \frac{1}{4} P(V_1=0 ,V_2=0,V_3=0,V_{r,1}=0, V_{r,2}=0, V_{r,3}=1|C_r = b, C_{n-r} = a)\\
&=& \frac{1}{4} g_{2,1, {\bf n}(b, a)}^{(r)}.
\end{eqnarray*}
Under $q=1$, cases 3 and 4 have probability zero, so
summing and applying Corollary \ref{cor:spectral.4} yields
\bea \nn
E(\xi_{n,3,1}^2)=E\left( \xi_{n,3,1,1}^2 + \xi_{n,3,1,2}^2 \right) = \frac{(n)_3}{4n^2(m+1)^2} \sum_{r = m}^{m+1} \sum_{a, b \in \{0, 1\}} \left( g_{1, 2, {\bf n}(b, a)}^{(r)} +  g_{2, 1, {\bf n}(b, a)}^{(r)}\right)\nn \\
= \frac{(n)_3}{2n^2(m+1)^2}\sum_{a, b \in \{0, 1\}} \left(g_{1, 2,{\bf n}(b, a)}^{(m)}+g_{2, 1, {\bf n}(b, a)}^{(m)}\right)\nn \\
=\frac{m(2m-1)}{4n^2(m+1)}\left( 1 - \overline{\lambda}_{n,2,{\bf n}_m} \right).\label{xi-n-1:three-distinct-case1}
\ena

Moving to $q=2$, for $E\xi_{n,3,2,1}^2$ we compute
\begin{eqnarray*}
\lefteqn{P(V_i=0 ,V_{j_1}=0,V_{j_2}=0,V_{r,i}=0, V_{t, i} = 0, V_{r,j_1}=1, V_{t,j_2}=1, C_r = b, C_{t} = a)}\\
&=& \frac{1}{4}P(V_1=0 ,V_2=0,V_3=0,V_{r,1}=0, V_{r,2}=1, V_{t,1} = 0, V_{t,3}=1|C_r = b, C_t = a) \\
&=& \frac{1}{4} g_{[(1, 0), (2, 1); (1, 0), (3, 1)], {\bf
n}(b, a)}^{(r,t)},
\end{eqnarray*}
while similarly, with factors of $1/4$ for $E\xi_{n,3,2,2}^2$ we have
\begin{eqnarray*}
P(V_1=0 ,V_2=0,V_3=0,V_{r,1}=1,V_{r,2}=0, V_{t,1}=1, V_{t,3}=0 |C_r = b, C_t = a)
=g_{[(1, 1), (2, 0); (1, 1), (3, 0)],{\bf n}(b, a)}^{(r,t)},
\end{eqnarray*}
for $E\xi_{n,3,2,3}^2$ we have
\beas
P(V_1=0 ,V_2=0,V_3=0, V_{r,1}=0, V_{r,2}=1, V_{t,1}=1 , V_{t,3}=0|C_r = b, C_t = a)=g_{[(1, 0), (2, 1); (1, 1), (3, 0)],{\bf n}(b, a)}^{(r,t)},
\enas
and for $E\xi_{n,3,2,4}^2$ we consider
\begin{eqnarray*}
P(V_1=0 ,V_2=0, V_3=0, V_{r,1}=1, V_{r,2}=0, V_{t,1}=0, V_{t,3}=1|C_r = b, C_t = a)= \frac{1}{4}g_{[(1, 1), (2, 0); (1, 0), (3, 1)], {\bf n}(b, a)}^{(r,t)}.
\end{eqnarray*}

Summing up these $q=2$ case terms and applying (\ref{last.4.xi-n-3:exp}) in the Appendix, we obtain
\begin{multline} \label{xi-n-3:exp}
E(\xi_{n,3,2}^2)=E(\xi_{n,3,2,1}^2+\xi_{n,3,2,2}^2+\xi_{n,3,2,3}^2+\xi_{n,3,2,4}^2) \\
= \frac{(n)_3}{4n^2(m+1)^2}\sum_{a, b \in \{0, 1\}} \sum_{r \ne t} \left(
g_{[(1, 0), (2, 1); (1, 0), (3, 1)], {\bf
n}(b, a)}^{(r,t)} +g_{[(1, 1), (2, 0); (1, 1), (3, 0)], {\bf n}(b, a)}^{(r,t)} \right.  \\
 + \left.  g_{[(1, 0), (2, 1); (1, 1), (3, 0)],{\bf n}(b, a)}^{(r,t)} + g_{[(1, 1), (2, 0); (1, 0), (3, 1)], {\bf n}(b, a)}^{(r,t)} \right)\\
 =\frac{(n)_3}{4n^2(m+1)^2} \left( \frac{4(m+1)_2^2(2m-1)^2}{(n)_3^2} \right) = \frac{m(2m-1)}{2n^3}.
\end{multline}

Now summing (\ref{xi-n-3:exp}) and (\ref{xi-n-1:three-distinct-case1}) we obtain
\bea
E\xi_{n,3}^2  = \frac{m(2m-1)}{4n^2(m+1)}\left(1- \overline{\lambda}_{n,2,{\bf
n}_m }\right) + \frac{m(2m-1)}{2n^3}. \label{Exin32.sq}
\ena

It remains to consider $E\xi_{n,2}^2$. Under the constraints $i_1 \ne j_1, i_2 \ne j_2$, sums over indices that satisfy $|\{i_1,i_2,j_1,j_2\}|=2$
can be partitioned into the following two cases:
\beas
1.\,\, i_1=i_2, j_1 = j_2 \qmq{and} 2.\,\, i_1=j_2, j_1 = i_2.
\enas

Consider first $q=1$. Under Case 1, by consideration of the term $P(V_i=0,V_j=0, V_{r,i}=0,V_{r,j}=1 ,C_r = b, C_t = a)$
and the application of Corollary \ref{cor:spectral.4}, we obtain
\bea \label{n211}
E\xi_{n,2,1,1}^2 = \frac{(n)_2}{4n^2(m+1)^2}\sum_{r=m}^{m+1} \sum_{a,b \in \{0,1\}} g_{1, 1,{\bf n}(b, a)}^{(r)}
=
\frac{m}{2n^2(m+1)} \left(1-\overline{\lambda}_{n,2,{\bf n}_m} \right)\frac{}{}.
\ena
Case 2 has zero probability.

Now let $q=2$. For stages $r \not = t$, under Case 1 we consider
\beas
P(V_i=0, V_j=0, V_{r,i}=0,V_{r,j}=1, V_{t,i}=0, V_{t,j}=1,C_r = b, C_t = a)
= \frac{1}{4}g_{[(1,0),(2,1);(1,0),(2,1)], {\bf n}(b, a)}^{(r,t)}.
\enas
Now, summing over the two ways one can achieve $r \not = t$, which yield identical terms, and using (\ref{g.n221}) from the Appendix,
we obtain
\bea
E\xi_{n,2,2,1}^2 &=& \frac{(n)_2}{2n^2(m+1)^2} \sum_{a,b \in \{0,1\}} g_{[(1,0),(2,1);(1,0),(2,1)], {\bf n}(b, a)}^{(r,t)} \nn \\
&=&\frac{(n)_2}{2n^2(m+1)^2}\left(1+2\lambda_{n,1,{\bf n}_{m,m+1}} + \lambda_{n,2,{\bf n}_{m,m+1} } \right)\frac{(m+1)_2^2}{(n)_2^2}\nn \\
&=&\frac{m}{4n^3}(1+2\lambda_{n,1,{\bf n}_{m,m+1}} + \lambda_{n,2,{\bf n}_{m,m+1} }). \label{n221}
\ena

Under Case 2 we consider
\beas
P(V_i=0, V_j=0, V_{r,i}=0,V_{r,j}=1, V_{t,i}=1 , V_{t,j}=0,C_r = b, C_t = a) = \frac{1}{4}g_{[(1,0),(2,1);(1,1),(2,0)], {\bf n}(b, a)}^{(r,t)}
\enas
and similarly arrive at
\bea \label{n222}
E\xi_{n,2,2,2}^2 = \frac{m}{4n^3}(1-2\lambda_{n,1,{\bf n}_{m,m+1}} + \lambda_{n,2,{\bf n}_{m,m+1} }).
\ena

Summing up the terms (\ref{n211}), (\ref{n221}) and (\ref{n222}) yields
\bea
E\xi_{n,2}^2 =  \frac{m}{2n^2(m+1)}(1 - \overline{\lambda}_{n,2,{\bf n}_m})+ \frac{m}{2n^3}(1+\lambda_{n,2,{\bf n}_{m,m+1} }).\label{Exin22.sq}
\ena

Finally, recalling
\beas
B=\mbox{Var}(\xi_n)= E\left( \xi_{n,4}^2+\xi_{n,3}^2+\xi_{n,2}^2 \right) -(E\xi_n)^2,
\enas
collecting terms (\ref{mean.xi.n}), (\ref{Exin42.sq}), (\ref{Exin32.sq}), and (\ref{Exin22.sq})  and simplifying we obtain
\beas
B
&=& \frac{m(8m +3)}{8n^3} +\frac{m(m^2+m-1)}{4n^2(m+1)}\overline{\lambda}_{n,2,{\bf n}_m} + \frac{(m+1)_2}{4n^3}\lambda_{n,2,{\bf n}_{m,m+1} } \\
& & +\frac{m^2(m-1)}{8n^2(m+1)}\overline{\lambda}_{n,4,{\bf n}_{m}} + \frac{(m)_2^2}{2n^2(n)_4}\lambda_{n,4,{\bf n}_{m,m+1}} - \frac{m^2}{4n^2} \overline{\lambda}_{n,2,{\bf n}_m}^2. 
\enas


Applying Lemma \ref{lam24bds}, we find for $n \ge 7$
\beas
B \le \frac{m(8m +3)}{8n^3} + \left(\frac{1}{16} +
\frac{1}{16n} + \frac{1}{64}+\frac{1}{64n^2} \right)e^{-n} \le \frac{m(8m +3)}{8n^3}+\frac{e^{-n}}{11}.
\enas

Now, using (\ref{final.A}),
\beas
\lefteqn{\sqrt{A \cdot B}} \\
&\le& \sqrt{\frac{3}{8n^3}\frac{m(8m +3)}{8n^3}} +
\sqrt{\frac{m(8m +3)}{8n^3}\frac{e^{-n}}{8n^2}}
+ \sqrt{ \frac{3}{8n^3} \frac{e^{-n}}{11}} +
\sqrt{\frac{e^{-n}}{8n^2} \frac{e^{-n}}{11} } \nn \\
&\le& \frac{5m+1}{8n^3} + \sqrt{\frac{1}{4n}
\frac{e^{-n}}{8n^2}}+\sqrt{\frac{3}{8n^3}\frac{e^{-n}}{11}}
+ \sqrt{\frac{e^{-2n}}{88n^2}} \nn \\
& = &\frac{5m+1}{8n^3}   +
\sqrt{\frac{1}{32n^3}}e^{-n/2} +
\sqrt{\frac{3}{88n^3}}e^{-n/2}+ \frac{e^{-n}}{8n} \nn \\
&\le& \frac{5m+1}{8n^3} + \frac{e^{-n/2}}{2n^{3/2}} +
\frac{e^{-n}}{8n}. 
\enas

Hence, from (\ref{improved.CS.bound}), again using $n \ge 7$,
\beas
\lefteqn{\mbox{Var}(E(V^s-V|{\cal F})) \le A+4B+4\sqrt{AB}}\\
&\le& \left( \frac{3}{8n^3} + \frac{e^{-n}}{8n^2} \right)+4\left( \frac{m(8m +3)}{8n^3}+\frac{e^{-n}}{11}\right) + 4 \left(\frac{5m+1}{8n^3} + \frac{e^{-n/2}}{2n^{3/2}} +
\frac{e^{-n}}{8n} \right)\\
&\le& \frac{1}{n}+ \frac{e^{-n/2}}{8}.
\enas
Taking square root and using $\sqrt{a+b} \le \sqrt{a} + \sqrt{b}$ now yields the upper bound (\ref{upDelta1}), thus completing the proof. \bbox

We now provide a bound for the normal approximation of $X$ in the odd case. We remark that using $V$, fewer error terms, and therefore a smaller bound, results when standardizing $X$ as in Theorem \ref{thm:main}, that is, not by its own mean and variance but by the (exponentially close) mean and variance of the closely coupled $V$.

\noindent {\em Proof of Theorem \ref{thm:main}: odd $n$.}
Letting $W=(X-n/2)/\sigma_V$ and $W_V=(V-n/2)/\sigma_V$, recalling
$|X-V| \le 2$ from (\ref{difVX.is.2}), we have
\bea \label{WWVdiff1}
|W-W_V| = \left|X-V \right|/\sigma_V \le 2/\sigma_V.
\ena
With $\Phi(z)=P(Z \le z)$ and
$$
C_V= \frac{n}{2\sigma_V^2}{\overline \Delta}_1 + 1.64 \frac{n}{\sigma_V^3}+
\frac{2}{\sigma_V},
$$
by (\ref{WWVdiff1}) and Theorem \ref{thm:normV} we obtain
\beas
P(W \le z)-\Phi(z)  &\le& P(W_V-2/\sigma_V \le z) - \Phi(z)\\
&=& P(W_V \le z+2/\sigma_V) - \Phi(z+2/\sigma_V)+\Phi(z+2/\sigma_V)-\Phi(z)\\
&\le& C_V + 2/(\sigma_V \sqrt{2 \pi}).
\enas
As a corresponding lower bound can be similarly demonstrated, the claim is shown. \bbox

\section{Spectral Decomposition}
\label{sec:spe}

In [\ref{zl}] the lightbulb chain was analyzed as a composition chain of multinomial type. Such chains in general are based on a $d \times d$ Markov transition matrix $P$ that describes the transition of a single particle in a system of $n$ identical particles, a subset of which is selected uniformly to undergo transition at each time step according to $P$.

In the case of the lightbulb chain there are $d=2$ states and the transition matrix $P$ of a single bulb is given by
\beas
P=\left[
\begin{array}{cc}
0 & 1 \\
1 & 0
\end{array}
\right],
\enas
where we let $e_0=(1,0)'$ and $e_1=(0,1)'$ denote the $0$ and $1$ states of the bulb, for off and on, respectively.
With $b \in \{0,1,\ldots,n\}$ let $P_{n,b,s}$ be the $2^b \times 2^b$ transition matrix of a subset of size $b$
of the $n$ total lightbulbs when $s$ of the $n$ bulbs are selected uniformly to be switched. Letting $P_{n,0,s}=1$ for all $n$ and $s$, and $I_2$ the $2 \times 2$ identity matrix, for $n \ge 1$ the matrix $P_{n,b,s}$ is given recursively by
\beas
P_{n,b,s}= \frac{s}{n} \left(P \otimes P_{n-1,b-1,s-1} \right)+ (1-\frac{s}{n})\left(I_2 \otimes P_{n-1,b-1,s}\right)
\quad \mbox{for $b \in \{1,\ldots,n\}$,}
\enas
as any particular bulb among the $b$ in the subset considered is selected with probability $s/n$ to undergo transition according to $P$, leaving the $s-1$ remaining switches to be distributed over the remaining $b-1$ of $n-1$ bulbs, and with probability $1-s/n$ the bulb is left unchanged, leaving
all the $s$ switches to be distributed.

The transition matrix $P$ is easily diagonalizable by the orthogonal matrix $T$ as
\bea \label{def:Gamma}
P=T'\Gamma T \qmq{where}
\quad
T=\frac{1}{\sqrt{2}}\left[
\begin{array}{cc}
1 & 1 \\
-1 & 1
\end{array}
\right] \quad \mbox{and} \quad
\Gamma=\left[
\begin{array}{cc}
1 & 0 \\
0 & -1
\end{array}
\right],
\ena
hence $P_{n,b,s}$ is diagonalized by
\bea \label{P=ViGV}
P_{n,b,s}= \otimes^b T' \Gamma_{n,b,s} \otimes^b T
\ena
where $\Gamma_{n,0,s}=1$, and is given the recursion
\bea \label{gnbsrec}
\Gamma_{n,b,s}= \frac{s}{n} \left(\Gamma \otimes \Gamma_{n-1,b-1,s-1} \right)+ (1-\frac{s}{n})\left(I_2 \otimes \Gamma_{n-1,b-1,s}\right)
\quad \mbox{for $b \in \{1,\ldots,n\}$.}
\ena
The next result describes the diagonal matrices $\Gamma_{n,b,s}$ more explicitly in terms of a sequence of vectors ${\bf a}_b$ of length $2^b$ for all $b \ge 1$ defined through the recursion
\bea \label{arecursion}
{\bf a}_b=({\bf a}_{b-1},{\bf a}_{b-1}+{\bf 1}_{b-1}) \qmq{for $b \ge 2$, with ${\bf a}_1=(0,1)$,}
\ena
where ${\bf 1}_b=(1,\ldots,1)$ is of size $2^b$. For example,
\bea \label{bfa.2.3}
{\bf a}_1=(0,1), \quad {\bf a}_2=(0,1,1,2) \qmq{and} {\bf a}_3=(0,1,1,2,1,2,2,3).
\ena
Letting $a_n$ be the $n^{th}$ term of the vector ${\bf a}_b$ for any $b$ satisfying $2^b \ge n$ results
in a well defined sequence $a_1,a_2,\ldots$.

\begin{lemma} 
For $n\in \{0,1,\ldots,\},b,s \in \{0,\ldots,n\}$, and $\lambda_{n,b,s}$ given by (\ref{def-lam}), the
matrix $\Gamma_{n,b,s}$ in (\ref{P=ViGV}) satisfies
\beas
\Gamma_{n,b,s}=\mbox{diag}(\lambda_{n,a_1,s},\ldots,\lambda_{n,a_{2^b},s}).
\enas
In particular, with ${\bf 0}_{2^{b-1}}$ the vector of zeros of length $2^{b-1}$, for $b \ge 1$
\beas
\Gamma_{n,b,s}=\mbox{diag}(\lambda_{n,a_1,s},\ldots,\lambda_{n,a_{2^{b-1}},s},{\bf 0}_{2^{b-1}})+\mbox{diag}({\bf 0}_{2^{b-1}},\lambda_{n,a_1+1,s},\ldots,\lambda_{n,a_{2^{b-1}}+1,s}).
\enas
\end{lemma}

For instance, from (\ref{bfa.2.3}), for $b=2$ we have
\beas
\Gamma_{n,2,s}=\mbox{diag}(\lambda_{n,0,s},\lambda_{n,1,s},\lambda_{n,1,s},\lambda_{n,2,s}),
\enas
and for $b=3$,
\beas
\Gamma_{n,3,s}=\mbox{diag}(\lambda_{n,0,s},\lambda_{n,1,s},\lambda_{n,1,s},\lambda_{n,2,s},\lambda_{n,1,s},\lambda_{n,2,s},
\lambda_{n,2,s},\lambda_{n,3,s}).
\enas

\noindent \proof As $a_1=0$ we have $\Gamma_{n,0,s}=1=\lambda_{n,0,s}$, so the lemma is true for $b=0$. For the inductive step, assuming the lemma is true for $b-1$, by (\ref{gnbsrec}) and the definition of $\Gamma$ from (\ref{def:Gamma}), it suffices to verify
\beas
\frac{s}{n}\lambda_{n-1,a,s-1}+(1-\frac{s}{n})\lambda_{n-1,a,s} &=& \lambda_{n,a,s} \quad \mbox{and} \\
-\frac{s}{n}\lambda_{n-1,a_,s-1}+(1-\frac{s}{n})\lambda_{n-1,a_,s} &=& \lambda_{n,a+1,s},
\enas
for all $a=0,1,\ldots$. To prove the first equality, by (\ref{def-lam}) we have
\beas
\frac{s}{n}\lambda_{n-1,a,s-1}+(1-\frac{s}{n})\lambda_{n-1,a,s} &=&\frac{s}{n}\sum_{t=0}^a {a \choose t} (-2)^t \frac{(s-1)_t}{(n-1)_t}+(1-\frac{s}{n})\sum_{t=0}^a {a \choose t} (-2)^t \frac{(s)_t}{(n-1)_t}\\
&=&\sum_{t=0}^a {a \choose t} (-2)^t \left( \frac{s}{n}\frac{(s-1)_t}{(n-1)_t}+\left(1-\frac{s}{n}\right)\frac{(s)_t}{(n-1)_t}\right)\\
&=&\sum_{t=0}^a {a \choose t} (-2)^t \left( \frac{(s)_{t+1}+(n-s)(s)_t}{(n)_{t+1}}\right)\\
&=&\sum_{t=0}^a {a \choose t} (-2)^t \left( \frac{(s)_t(s-t+n-s)}{(n)_{t+1}}\right)\\
&=&\sum_{t=0}^a {a \choose t} (-2)^t \left( \frac{(s)_t(n-t)}{(n)_{t+1}}\right)\\
&=&\sum_{t=0}^a {a \choose t} (-2)^t \frac{(s)_t}{(n)_t}\\
&=&\lambda_{n,a,s}.
\enas
The second equality can be shown in similar, though slightly more involved, fashion.
\ignore{Proof of Second Equality:
\beas
\lefteqn{-\frac{s}{n}\sum_{t=0}^a {a \choose t} (-2)^t \frac{(s-1)_t}{(n-1)_t}+(1-\frac{s}{n})\sum_{t=0}^a {a \choose t} (-2)^t \frac{(s)_t}{(n-1)_t}}\\
&=&\sum_{t=0}^a {a \choose t} (-2)^t \left( -\frac{s}{n}\frac{(s-1)_t}{(n-1)_t}+(1-\frac{s}{n})\frac{(s)_t}{(n-1)_t}\right)\\
&=&\sum_{t=0}^a {a \choose t} (-2)^t \left( \frac{-(s)_{t+1}+(n-s)(s)_t}{(n)_{t+1}}\right)\\
&=&\sum_{t=0}^a {a \choose t} (-2)^t \left( \frac{(s)_t(t-s+n-s)}{(n)_{t+1}}\right)\\
&=&\sum_{t=0}^a {a \choose t} (-2)^t \left( \frac{(s)_t(n+t-2s)}{(n)_{t+1}}\right)\\
&=&\sum_{t=0}^a {a \choose t} (-2)^t \left( \frac{(s)_t(n-t-2(s-t))}{(n)_{t+1}}\right)\\
&=&\sum_{t=0}^a {a \choose t} (-2)^t \frac{(s)_t}{(n)_t} +\sum_{t=0}^a {a \choose t} (-2)^{t+1} \frac{(s)_{t+1}}{(n)_{t+1}}\\
&=&\sum_{t=0}^a {a \choose t} (-2)^t \frac{(s)_t}{(n)_t} +\sum_{t=1}^{a+1} {a \choose t-1} (-2)^t \frac{(s)_t}{(n)_t}\\
&=&\sum_{t=0}^{a+1} {a+1 \choose t} (-2)^t \frac{(s)_t}{(n)_t}\\
&=&\lambda_{n,a+1,s}.
\enas
}
\bbox

We note that [\ref{zl}] expresses these eigenvalues in terms of the hypergeometric function.

If the $k$ stages of the process $1,\ldots,k$ use switches ${\bf
s}=(s_1,\ldots,s_k)$, then
since the matrices $P_{n,b,s},s \in \{0,1,\ldots,n\}$ are simultaneously diagonalizable by
 (\ref{P=ViGV}),
the transition matrix $P_{n,b,{\bf s}}$ for any subset of $b$
bulbs can be diagonalized as
\bea
\label{direct}
P_{n,b,{\bf s}}=\prod_{j=1}^k P_{n,b,s_j}
=\otimes^b T' \Gamma_{n,b,{\bf s}} \otimes^b T = \otimes^b
T' \mbox{diag}(\lambda_{n,a_1,{\bf
s}},\ldots,\lambda_{n,a_{2^b},{\bf s}}) \otimes^b T,
\ena
where $\lambda_{n,a,{\bf s}}$ is given in (\ref{def-lam}) and
\bea
\label{deflambdabfs}
\Gamma_{n,b,{\bf s}}=\prod_{j=1}^k
\Gamma_{n,b,s_j}.
\ena
If  $\pi$ is a permutation of
$\{1,\ldots,k\}$ let $\pi ({\bf s})=
(s_{\pi(1)},\ldots,s_{\pi(k)})$. As all matrices $\Gamma_{n,b,s}$ are diagonal, from (\ref{deflambdabfs}) we have
$\Gamma_{n,b,{\bf s}}=\Gamma_{n,b,\pi({\bf s})}$, and now from
(\ref{direct}) we have the following result.
\begin{lemma}
\label{lem:permuation}
The distribution of the lightbulb chain is independent of the order in which the switch
variables ${\bf s}$ are applied, that is, for all permutations $\pi$,
\beas 
P_{n,b,{\bf s}}=P_{n,b,\pi({\bf
s})}.
\enas
\end{lemma}

The following lemma helps us compute probabilities such as $g_{\alpha,\beta,{\bf s}}^{(l)}$ in
(\ref{def:galphabeta}). For $j \in \{0,1,\ldots\}$ let $\Omega_{b,j}$ be the $2^b \times
2^b$ diagonal matrix in the variables $x_k, k \in \{0,1,\ldots\}$
given by
\bea \label{def:Omegabj}
\Omega_{b,j}=\mbox{diag}(x_{a_1+j},\ldots,x_{a_{2^b}+j}), \qmq{and set}
u_b=\otimes^b T e_0^{\otimes b} \qmq{and} w_b= \otimes^b T e_1^{\otimes b},
\ena
where we recall $e_0=(1,0)'$ and $e_1=(0,1)'$. Note that for $b=1$ we have
\bea \label{u1w1}
u_1=Te_0 = \frac{1}{\sqrt{2}}(1,-1)' \qmq{and} w_1=Te_1 = \frac{1}{\sqrt{2}}(1,1)'.
\ena

\begin{lemma} \label{lem:om.ind}
Let $t \in \{0,1,\ldots\}$ and $\Omega_t=\Omega_{t,0}$, and suppose that for some vector $v_t \in \mathbb{R}^{2^t}$, that
\bea
\label{wOr-1.0}
v_t' \Omega_t u_t = \frac{1}{2^t} \sum_{j=0}^t  {t
\choose j} a(j) x_j
\ena
holds for $t=b-1$ with some sequence $a(j),j=0,\ldots,b-1$. Then (\ref{wOr-1.0}) holds for $t=b$ when replacing $v_t$ by
$v_b=u_1 \otimes v_t$ and $a(j)$ by
\bea \label{b+}
a_u(j)=\left( \frac{b-j}{b} \right)a(j) + \left( \frac{j}{b} \right)a(j-1),
\ena
and for $t=b$ when replacing $v_t$ by $v_b=w_1 \otimes v_t$ and $a(j)$ by
\bea \label{b-}
a_w(j)=\left( \frac{b-j}{b} \right)a(j) - \left( \frac{j}{b} \right)a(j-1).
\ena

\end{lemma}

\noindent \proof By (\ref{arecursion}) we may write
\beas 
\Omega_b=\left[
\begin{array}{cc}
\Omega_{b-1,0} & 0 \\
0 & \Omega_{b-1,1}
\end{array}
\right]
\enas
and by (\ref{u1w1}) we have
\beas 
u_b = u_1 \otimes u_{b-1}= \frac{1}{\sqrt{2}}(u_{b-1}',-u_{b-1}')'.
\enas
Hence, when $v_b=u_1 \otimes v_{b-1}=(v_{b-1}',-v_{b-1}')'/\sqrt{2}$, we obtain
\beas
v_b'\Omega_b u_b&=& \frac{1}{2} \left( v_{b-1}'\Omega_{b-1,0}u_{b-1} + v_{b-1}'\Omega_{b-1,1}u_{b-1} \right) \\
&=& \frac{1}{2^b}\sum_{j=0}^{b-1}{b-1 \choose j}a(j)x_j+\frac{1}{2^b}\sum_{j=0}^{b-1}{b-1 \choose j}a(j)x_{j+1}\\
&=& \frac{1}{2^b}\sum_{j=0}^{b-1}{b-1 \choose j}a(j)x_j+\frac{1}{2^b}\sum_{j=1}^b{b-1 \choose j-1}a(j-1)x_j\\
&=& \frac{1}{2^b}\sum_{j=0}^b \left( {b-1 \choose j}a(j)+ {b-1 \choose j-1}a(j-1) \right)x_j\\
&=& \frac{1}{2^b}\sum_{j=0}^b {b \choose j}\left( \left(\frac{b-j}{b}\right)a(j) + \left(\frac{j}{b}\right)a(j-1) \right)x_j\\
&=& \frac{1}{2^b}\sum_{j=0}^b {b \choose j}a_u(j)x_j,
\enas
as claimed. The proof is essentially the same, using (\ref{u1w1}), when $v_b=w_1 \otimes v_{b-1} =(v_{b-1}',v_{b-1}')'/\sqrt{2}$.
\bbox

In Corollary \ref{cor:spectral.4} we express the functions $g_{\alpha,\beta,{\bf s}}^{(l)}$ in (\ref{def:galphabeta}) using the corresponding conditional probabilities
\bea \label{def:falphabeta}
f_{\alpha,\beta,{\bf s}}=P(X_i=0, i=1,\ldots,\alpha+\beta |X_{0,i} = 0, i=1,\ldots, \alpha ,X_{0,i}=1, i=\alpha+1,\ldots,\alpha+\beta)
\ena
that are the subject of the next lemma. As with $g_{\alpha,\beta,{\bf s}}^{(l)}$, we drop the dependence on ${\bf s}$ when ${\bf s}={\bf n}$.

\begin{lemma}
\label{lem:aor.spec} For given $\alpha,\beta \ge 0$, setting $b=\alpha+\beta$ the probability $f_{\alpha,\beta,{\bf s}}$ in (\ref{def:falphabeta}) is given by
\bea \label{falphabeta.prob.sum}
f_{\alpha,\beta,{\bf s}}=
\frac{1}{2^b}\sum_{j=0}^b {b \choose j}a_{\alpha,\beta}(j) \lambda_{n,j, {\bf s}},
\ena
where $a_{\alpha,0}(j)=1$ for all $\alpha \ge 0$ and
\bea \label{lem:aor.spec.induction}
a_{\alpha,\beta}(j)=\left( \frac{b-j}{b} \right)a_{\alpha,\beta-1}(j) -
\left( \frac{j}{b} \right)a_{\alpha,\beta-1}(j-1)
\quad \mbox{for all $\alpha \ge 0, \beta \ge 1$.}
\ena
\end{lemma}

\noindent \proof
Using exchangeability for the first equality, extracting the relevant component of the $k$-step transition matrix and applying (\ref{direct}) we obtain
\beas
f_{\alpha,\beta,{\bf s}}&=&
\left( (e_1')^{\otimes \beta} \otimes (e_0')^{\otimes \alpha} \right)
P_{n,b,{\bf s}}e_0^{\otimes b}
= \left( (e_1')^{\otimes \beta} \otimes (e_0')^{\otimes \alpha} \right)\otimes^b T' \Gamma_{n,b,{\bf s}} \otimes^b T e_0^{\otimes b}=v_b'\Gamma_{n,b,{\bf s}} u_b,
\enas
where $u_\alpha$ and $w_\beta$ are given as in (\ref{def:Omegabj}) and $v_b =w_\beta \otimes u_\alpha$.
Hence, with $\Omega_b=\Omega_{b,0}$ as in (\ref{def:Omegabj}), the result follows from
\bea \label{genvrur}
v_b'\Omega_b u_b = \frac{1}{2^b}\sum_{j=0}^b {b \choose j}a_{\alpha,\beta}(j)x_j.
\ena

We first prove the case $\beta=0$ by induction in $\alpha$; note in this case $v_b=u_b$. Equality (\ref{genvrur}) holds
with $a_{\alpha,0}(j)=1$ for $\alpha=0$, as both sides equal $x_0$ in this case. Assuming that (\ref{genvrur}) holds for
some $\alpha \ge 0$ with $a_{\alpha,0}(j)=1$,
then (\ref{b+}) of Lemma \ref{lem:om.ind} implies that
(\ref{genvrur}) holds for $\alpha+1$ and $\beta=0$ with
\beas
a_{\alpha+1,0}(j) = \left( \frac{b-j}{b} \right) a_{\alpha,0}(j) +  \left( \frac{j}{b} \right) a_{\alpha,0}(j-1)=1.
\enas
Hence (\ref{genvrur}) holds for all $\alpha \ge 0$ and $\beta=0$ with $a_{\alpha,0}(j)=1$.
Similarly, assuming now that (\ref{genvrur}) holds for $a_{\alpha,\beta-1}(j)$ with nonnegative $\alpha,\beta-1$, we have that (\ref{genvrur}) holds
with $a_{\alpha,\beta}(j)$ given by (\ref{lem:aor.spec.induction}) by (\ref{b-}) of Lemma \ref{lem:om.ind}, thus completing the induction.
\bbox

As our computations involve only moments up to fourth order, we highlight these particular special cases of
Lemma \ref{lem:aor.spec} in the following Corollary.
\begin{corollary} \label{cor:spectral.4}
For $\alpha,\beta \ge 0$, the probability $g_{\alpha,\beta,{\bf s}}^{(l)}$ in (\ref{def:galphabeta}) is given by
\bea \label{probtau2}
g_{\alpha,\beta,{\bf s}}^{(l)}=f_{\alpha,\beta,{\bf s}_l}\,p_{\alpha,\beta}({\bf s},l) \qmq{where} p_{\alpha,\beta}({\bf s},l) = \frac{(n-s_l)_\alpha (s_l)_\beta}{(n)_{\alpha+\beta}},
\ena
and
$f_{\alpha,\beta,{\bf s}}$ is given by (\ref{falphabeta.prob.sum}).  For $0 \le \alpha + \beta \le 4$,
the sequences in (\ref{falphabeta.prob.sum}) specialize to
\beas
a_{0,0}(j)=1
\enas
\beas
a_{0,1}(j)=(-1)^j \qmq{and} a_{1,0}(j)=1
\enas
\beas
a_{0,2}(j)=(-1)^j,\quad a_{1,1}(j)= 1-j \qmq{and} a_{2,0}(j)=1
\enas
\beas
a_{0,3}(j)=(-1)^j,\quad a_{1,2}(j)= (-1)^j(1-2j/3), \quad a_{2,1}=1-2j/3 \qmq{and} a_{3,0}(j)=1
\enas
and
\beas
a_{0,4}(j)=(-1)^j,\quad a_{1,3}(j)= (-1)^j(1-j/2), \quad a_{2,2}=1-j(4-j)/3, \quad a_{3,1}=1-j/2 \qmq{and} a_{4,0}(j)=1.
\enas
\end{corollary}

\proof By Lemma \ref{lem:permuation}, that is, the fact that the switch variables can be
applied in any order, conditioning on the values of the switch variables
in stage $l$ yields the same result as assuming these values as initial
conditions in stage 0, and applying the switch pattern ${\bf s}_l$,
that is, ${\bf s}$ skipping stage $l$. Hence the first claim in (\ref{probtau2}) follows, as
the first factor is the probability of the given event conditioned on the values in stage $l$,
while the second factor is the probability of the conditioning event, as
\beas
P(X_{l1}=\cdots=X_{l\alpha}=0,X_{l,\alpha+1}=\cdots=X_{l,\alpha+\beta}=1)
=\prod_{i=0}^{\alpha-1} \left( \frac{n-s_l-i}{n-i}\right)   \prod_{i=0}^{\beta-1} \left( \frac{s_l-i}{n-\alpha-i} \right) =p_{\alpha,\beta}({\bf s},l).
\enas

The specific forms of the sequences $a_{\alpha,\beta}(j)$ for $0 \le \alpha + \beta \le 4$ follow directly from the initial condition and recursion in Lemma \ref{lem:aor.spec}.
\bbox

Applying Corollary \ref{cor:spectral.4}, we obtain, for example, the formulas
\beas 
g_{2,1,{\bf s}}^{(l)}= \frac{1}{8}(1+\lambda_{n,1,{\bf s}_l }-\lambda_{n,2,{\bf s}_l }-\lambda_{n,3,{\bf s}_l }) \frac{(n-s_l)_2 s_l}{(n)_3},
\enas
and
\beas 
g_{2,2,{\bf s}}^{(l)}= \frac{1}{16}(1-2\lambda_{n,2,{\bf s}_l }+\lambda_{n,4,{\bf s}_l }) \frac{(s_l)_2(n-s_l)_2}{(n)_4}.
\enas

Lastly we present the bounds on products of eigenvalues of the chain used to handle the variance term
(\ref{def-Delta}) when applying Theorem \ref{thm:main} to the lightbulb chain.
\begin{lemma}
\label{lam24bds} For all even $n \ge 6$, and ${\bf s} = {\bf
n}_{n/2}$, or for all
odd $n=2m+1 \ge 7$ and ${\bf s} \in \{{\bf n}_m,{\bf n}_{m+1}, {\bf
n}_{m,m+1} \}$,
\beas 
|\lambda_{n, 2, {\bf s}}| \le e^{-n}  \qmq{and} |\lambda_{n, 4,{\bf
s}}| \le \frac{1}{2}e^{-n}.
\enas
Moreover, with $\overline{\lambda}_{n,b,{\bf s}}$ as given in (\ref{def:overlinelambda}),
for all odd $n =2m+1 \ge 7$ we have
\beas 
|\overline{\lambda}_{n, 2, {\bf n}_m}| \le e^{-n}  \qmq{and}
|\overline{\lambda}_{n, 4,{\bf n}_m}| \le \frac{1}{2}e^{-n}.
\enas
\end{lemma}

\noindent \proof The following calculations slightly generalize the arguments of [\ref{rrz}].
Let $m$ be given by $n=2m$, and $n=2m+1$, in the even and odd cases, respectively.
For $n \ge 2$ consider the second degree polynomial (cf. (\ref{lambdan12s}))
\beas
f_2(x)= 1-\frac{4x}{n}+\frac{4(x)_2}{(n)_2}, \quad 0 \le x \le n.
\enas
It is simple
to verify that $f_2(x)$ achieves its global minimum value of $-1/(n-1)$ at $n/2$,
and that $f_2(x)$ has exactly two roots, at $(n +
\sqrt{n})/2$ and $(n - \sqrt{n})/2$. Hence, as $f_2(x) \le 0$ for all $x$ between these roots, and
additionally, as $(x-1)/(n-1) \le x/n$ for all $x \in [0,n]$, we obtain the bound
\bea \label{lam2bounds}
\left|f_2(x)\right| \le
\left\{
\begin{array}{cl}
\frac{1}{n-1} & \mbox{for $x \in [\frac{n -
\sqrt{n}}{2},\frac{n +
\sqrt{n}}{2}]$} \\
\left(1 - \frac{2x}{n}\right)^2 & \mbox{for $x \in [\frac{n -
\sqrt{n}}{2},\frac{n +
\sqrt{n}}{2}]^c \cap [0,n]$.}
\end{array}
\right.
\ena

For $x \in \mathbb{R}$ let $\lfloor x \rfloor$ and $\lceil x \rceil$
denote the greatest integer less than or equal to $x$, and the
smallest integer greater than or equal to $x$. In both the even and odd cases let
\beas
{\bf t}= \left\{\lceil \frac{n-\sqrt{n}}{2} \rceil, \cdots,
\lfloor \frac{n+\sqrt{n}}{2} \rfloor \right\}\setminus \{m,m+1\}
\enas
so that
\beas
\vert \lambda_{n,2,{\bf n}_{m,m+1}} \vert &=& \left(
\prod_{s=0}^{\lfloor \frac{n-\sqrt{n}}{2} \rfloor}\vert f_2(s) \vert
\right) \left( \prod_{s \in {\bf t}} \vert f_2(s) \vert \right)
\left( \prod_{s=\lceil \frac{n+\sqrt{n}}{2}\rceil}^{n}\vert f_2(s)
\vert \right).
\enas
If either of the roots
$(n-\sqrt{n})/2$ or $(n+\sqrt{n})/2$ is an integer then equality
holds as both expressions above are zero. Now assuming
neither value is an integer, the product is over disjoint indices.

Applying the bound (\ref{lam2bounds}), $\lfloor n/2-x \rfloor + \lceil n/2+x \rceil = n$ and
$1-2(n-s)/n=-(1-2s/n)$
 yields
\beas
\vert \lambda_{n,2,{\bf n}_{m,m+1}} \vert &\le&
\left( \prod_{s=0}^{\lfloor \frac{n-\sqrt{n}}{2} \rfloor} \left(1-\frac{2s}{n} \right)^2  \right)
\left( \prod_{s \in {\bf t}}  \frac{1}{n-1} \right)
\left( \prod_{s=\lceil \frac{n+\sqrt{n}}{2}\rceil}^n \left(1-\frac{2s}{n} \right)^2 \right) \\
&=&
\left( \prod_{s=0}^{\lfloor \frac{n-\sqrt{n}}{2} \rfloor} \left(1-\frac{2s}{n} \right)  \right)^4
\left(\frac{1}{n-1} \right)^{|{\bf t}|},
\enas
where $|{\bf t}|$ is the cardinality of ${\bf t}$.

Using $1-x \le e^{-x}$ and that $\lfloor x \rfloor \ge x-1$ on the first product,
we obtain the bound
\beas 
\vert \lambda_{n,2,{\bf n}_{m,m+1}} \vert
&\le& \left( e^{-\frac{2}{n}(\lfloor{\frac{n-\sqrt{n}}{2}}\rfloor)(\lfloor{\frac{n-\sqrt{n}}{2}}\rfloor+1)/2} \right)^4
e^{-|{\bf t}| \log (n-1)} \le e^{-(n-2\sqrt{n}-1+2/\sqrt{n}+|{\bf t}| \log (n-1))}.
\enas

To control $|{\bf t}|$, note that as $\lceil x \rceil \le x+1$, we have
\beas
|{\bf t}|= \lfloor{\frac{n+\sqrt{n}}{2}}\rfloor-\lceil{\frac{n-\sqrt{n}}{2}}\rceil  -1 \ge \sqrt{n}-3.
\enas
As $\log 63 \ge 4$, the result follows for all $n \ge 64$ from
\beas
-2\sqrt{n}-1+2/\sqrt{n}+|{\bf t}|\log (n-1) &\ge&
-2\sqrt{n}-1+4(\sqrt{n}-3) = 2\sqrt{n}-13 \ge 0,
\enas
and by direct verification for all odd $7 \le n \le 63$, thus proving the claimed inequality
for $\lambda_{n,2,{\bf n}_{m,m+1}}$ for all odd $n \ge 7$.

By (\ref{lam2bounds}) we have $|f_2(x)|\le 1$ for $x \in [0,n]$, and hence
\bea \label{lamb2.to.lamb4}
\left| \lambda_{n,2,{\bf n}_m} \right| = \left| f_2(m+1) \lambda_{n,2,{\bf n}_{m,m+1}}\right| \le \left|\lambda_{n,2,{\bf n}_{m,m+1}}\right|,
\ena
so the claimed bound holds for $\lambda_{n,2,{\bf n}_m}$, and
likewise for $\lambda_{n,2,{\bf n}_{m+1}}$, for all odd $n \ge 7$.
The claim for $\overline{\lambda}_{n,2,{\bf n}_m}=(\lambda_{n,2,{\bf
n}_m}+\lambda_{n,2,{\bf n}_{m+1}})/2$ now follows. Direct
verification of the claim for $\lambda_{n,2,{\bf n}_m}$ for all
even $6 \le n \le 62$ now also completes the proof for all
cases involving $\lambda_{n,2,{\bf s}}$.

Now turning to $\lambda_{n, 4, {\bf s}}$, for all $n \ge 6$ consider the fourth degree polynomial
\begin{eqnarray*}
f_4(x) = 1-\frac{8x}{n}+\frac{24(x)_2}{(n)_2}-\frac{32(x)_3}{(n)_3}+\frac{16(x)_4}{(n)_4}, \quad \quad 0 \le x \le n.
\end{eqnarray*}
It can be checked that the four roots of $f_4(x)$ are given by
\beas
x_{1\pm} = \frac{n\pm\sqrt{\sqrt{2} \sqrt{3n^2-9n+8}+3n-4}}{2}
\qmq{and} x_{2\pm} =
\frac{n\pm\sqrt{-\sqrt{2} \sqrt{3n^2-9n+8}+3n-4}}{2},
\enas
and that additionally the three roots to the cubic equation $f_4'(x)=0$ occur at
\beas
y_1=n/2 \qmq{and} y_{2\pm}=(n\pm \sqrt{3n-4})/2.
\enas
These roots satisfy
\beas
0 < x_{1-} <y_{2-} < x_{2-} < y_1 <  x_{2+} < y_{2+}< x_{1+} < n.
\enas
To obtain a bound over the interval $[x_{1-},x_{1+}]$, evaluating $f_4(x)$ at its
critical values we obtain
\beas
f_4(y_1) &=& \frac{3}{(n-1)(n-3)} \le \frac{3}{(n-3)^2}\quad \mbox{and} \\
f_4(y_{2\pm}) &=& -\frac{2(3n^2-9n+8)}{n(n-1)(n-2)(n-3)} \ge
-\frac{6}{(n-3)^2}. \enas

To bound $f_4(x)$ by $f_2^2(x)$ in the remaining part of $[0,n]$, write
\beas
f_2^2(x)-f_4(x) = \frac{16 x (n-x) p(x)}{(n-1)^2 n^2 \left(n^2-5 n+6\right)}
\enas
where
\beas
p(x)=(4n-6)x^2 + (6n-4n^2)x+ n^3-2 n^2+n.
\enas
The roots of the quadratic $p(x)$ are given by
\beas
z_{\pm} = \frac{n}{2} \pm \frac{1}{2} \sqrt{\frac{n(n-2)}{2n-3}}.
\enas
As $5n^2-15n+12 \ge 0$ for all $n$, we have $(2n-3)(3n-4) \ge n(n-2)$,
and therefore
\beas
\left( \sqrt{2}\sqrt{3n^2-9n+8}+3n-4\right)(2n-3) \ge n(n-2).
\enas
Dividing by $2n-3$ and taking square roots demonstrates that
\beas
x_{1-} \le z_{-}<z_{+} < x_{1+}.
\enas
Hence $p(x)$ is nonnegative on the complement of $[z_{-},z_{+}]$, and we obtain
\bea \label{lam4bounds}
|f_4(x)| \le \left\{
\begin{array}{cl}
\frac{6}{(n-3)^2}& \mbox{for $x \in [x_{1-},x_{1+}]$} \\
f_2^2(x) & \mbox{for $x \not \in [x_{1-},x_{1+}], x \in [0,n]$.}
\end{array}
\right.
\ena

Now write for short
$$
C(n) = \sqrt{\sqrt{2} \sqrt{3n^2-9n+8} + 3n - 4} \qmq{so that} x_{1-} = \frac{n - C(n)}{2}.
$$
Using (\ref{lam4bounds}), (\ref{lam2bounds}) and that $|s-n/2| > \sqrt{n}/2$ whenever $s \le (n-C(n))/2$,
$1 - x \le e^{-x}$
for $x \ge 0$ and finally $\lfloor{x}\rfloor \ge x - 1$, we obtain
\bea \label{l4cnleft}
\prod_{s = 0}^{\lfloor{\frac{n - C(n)}{2}\rfloor}} |\lambda_{n, 4, s}| \le
\prod_{s = 0}^{\lfloor{\frac{n - C(n)}{2}\rfloor}}
f_2^2(s) \le  \left( \prod_{s =
0}^{\lfloor{\frac{n - C(n)}{2}} \rfloor}
\left(1-\frac{2s}{n}\right) \right)^4 \le
e^{-(n-2(C(n)+1) + C(n)^2/n + 2C(n)/n)}.
\ena
Now let
\beas
{\bf u}=
\left\{\lceil \frac{n-C(n)}{2} \rceil, \cdots ,
\lfloor \frac{n+C(n)}{2} \rfloor \right\}\setminus \{m,m+1\}.
\enas

Using $\lambda_{n,4,s}=\lambda_{n,4,n-s}$, (\ref{lam4bounds}), $|{\bf u}| \ge C(n) - 3$ and (\ref{l4cnleft}), we have that
\beas
\nn \lefteqn{|\lambda_{n, 4, {\bf n}_{m,m+1}}|}\\
&=& \prod_{s =
0}^{\lfloor{\frac{n - C(n)}{2}\rfloor}} |\lambda_{n, 4, s}| \prod_{s \in
{\bf u}} |\lambda_{n, 4, s}| \prod_{s =\lceil \frac{n+C(n)}{2} \rceil}^n
|\lambda_{n, 4, s}| \le \left(\prod_{s = 0}^{\lfloor{\frac{n - C(n)}{2}\rfloor}}
\lambda_{n, 4, s}^2 \right)\left(
\frac{6}{(n-3)^2}\right)^{C(n)-3}\\
&\le& e^{-2(n-2(C(n)+1)+ C(n)^2/n + 2C(n)/n)+(C(n)-3)(\log(6)-2\log (n-3))}\nn \\
\nn &=& e^{-n} \cdot e^{-(n - 4(C(n)+1)+ 2C(n)^2/n + 4C(n)/n +
(C(n)-3)(2\log(n-3)-\log(6)))}.
\enas
Since for $n \ge 96$ we have
\beas
n-4(C(n)+1) \ge \log(2) \qmq{and} (C(n)-3)\left( 2\log (n-3) -  \log(6)\right) \ge 0,
\enas
we conclude that $|\lambda_{n, 4, {\bf n}_{m,m+1}}| \le e^{-n}/2$ for all such $n \ge 96$, and direct
verification shows this same bound holds for all $7 \le n \le 95$.  Arguing as in
(\ref{lamb2.to.lamb4}), the claimed inequality therefore holds for
$\lambda_{n,4,{\bf n}_m}$ and $\lambda_{n,4,{\bf n}_{m+1}}$, and hence also for
$\overline{\lambda}_{n,4,{\bf n}_m}=(\lambda_{n,4,{\bf
n}_m}+\lambda_{n,4,{\bf n}_{m+1}})/2$ for all
odd $n \ge 7$, as well as for $\lambda_{n,4,{\bf n}_{n/2}}$ for all even $n \ge 8$.
The proof is completed by directly verifying the bound for $\lambda_{n,4,{\bf n}_{n/2}}$
for $n=6$.
\bbox

\section*{Appendix}
We work out some of the detailed calculations in the proof of Theorem \ref{thm:normV}. First, we show
\begin{eqnarray}
\lefteqn{\sum_{a,b \in \{0,1\}} g_{[(1,0), (2, 1); (3, 0), (4, 1)],
{\bf n}(b, a)}^{(m,m+1)}} \nonumber \\
&=& \left( \frac{(m+1)_3}{(n)_4}\right)^2 \left( (2m-1)^2 + 2(2m-1)\lambda_{n,2,{\bf n}_{m,m+1}} + \lambda_{n,4,{\bf n}_{m,m+1}} \right).
\label{g:xi-42}
\end{eqnarray}

To handle this sum, first write
\begin{multline} \label{xi-n-4:exp.sum}
\sum_{a,b \in \{0,1\}} g_{[(1,0), (2, 1); (3, 0), (4, 1)],
{\bf n}(b, a)}^{(m,m+1)}\\
= \sum_{a,b \in \{0,1\}} \sum_{e_1, e_2, e_3, e_4 \in \{0, 1\}} P\left(X_i=0, i=1,\ldots,4,X_{m1}=0, X_{m2} = 1, X_{m3}=e_1, X_{m4} = e_2 \right.\\
 \left. X_{m+1,1} = e_3,  X_{m+1, 2} = e_4 , X_{m+1,3} = 0, X_{m+1, 4} = 1\right).
\end{multline}
We may write the probability in the sum as the conditional probability of $\{X_i=0,i=1,\ldots,4\}$ given
$\{(X_{m1}, X_{m2}, X_{m3}, X_{m4}) = (0,1,e_1,e_2), (X_{m+1,1},  X_{m+1, 2}, X_{m+1,3}, X_{m+1, 4})= (e_3,e_4,0,1)\}$,
multiplied by the unconditional probability of this last event. By Lemma \ref{lem:permuation}, that is, that the switch variables may be applied in any order, conditioning the first event on the values of the switch variables for these four bulbs in stages $m$ and $m+1$ is equivalent to conditioning on the combined event $\{(X_{01},X_{02},X_{03},X_{04})=(e_3, e_4+1, e_1, e_2+1)\}$ in, say, an initial stage 0, and running the lightbulb process in the remaining $n-2$ stages, with, recalling our extension of the notation defined in (\ref{sdell}), pattern ${\bf n}_{m,m+1}=(1,\ldots,m-1,m+2,\ldots,n)$.

For the conditional probability, by (\ref{def:falphabeta}) we may write
\bea \label{use.for.3}
P\left(X_i=0, i=1,\ldots,4| (X_{01},X_{02},X_{03},X_{04})=(e_3, e_4+1,  e_1,   e_2+1)\right)  = f_{\alpha,\beta,{\bf n}_{m,m+1}}
\ena
where $\alpha+\beta=4$ and
\beas
\alpha={\bf 1}(e_3=0)+{\bf 1}(e_4+1=0)+{\bf 1}(e_1=0)+{\bf 1}(e_2+1=0).
\enas
Having conditioned on the switch variables in stages $m$ and $m+1$, we note that this conditional probability does not depend on $a$ or $b$.

For the unconditional probability, we have
\begin{multline*}
P \left((X_{m1}, X_{m2}, X_{m3}, X_{m4}) = (0,1,e_1,e_2), (X_{m+1,1},  X_{m+1, 2}, X_{m+1,3}, X_{m+1, 4})= (e_3,e_4,0,1)\right)\\
=
\left( \frac{(m+1-b)_{3-e_1-e_2}(m+b)_{1+e_1+e_2}}{(n)_4} \right) \left( \frac{(m+1-a)_{3-e_3-e_4}(m+a)_{1+e_3+e_4}}{(n)_4} \right),
\end{multline*}
as the stages are independent, and the event requires $1+e_1+e_2$ switches to be drawn from $m+b$ in stage $m$, and $1+e_3+e_4$ switches to be drawn from $m+a$ in stage $m+1$. Hence, changing the order of summation in (\ref{xi-n-4:exp.sum}), and recalling that the conditional probability $f_{\alpha,\beta,{\bf n}_{m,m+1}}$ does not depend on $a,b$, we have
\bea \label{sec.sum.xi-n-4:exp}
\sum_{a,b \in \{0,1\}} g_{[(1,0), (2, 1); (3, 0), (4, 1)],
{\bf n}(b, a)}^{(m,m+1)}=
\sum_{\alpha=0}^4 f_{\alpha,\beta,{\bf n}_{m,m+1}}w_{\alpha,\beta}
\ena
where
\beas
w_{\alpha,\beta}
=\sum_{E_{\alpha,\beta}}\sum_{\{a,b\}\in \{0,1\}} \left( \frac{(m+1-b)_{3-e_1-e_2}(m+b)_{1+e_1+e_2}}{(n)_4} \right) \left( \frac{(m+1-a)_{3-e_3-e_4}(m+a)_{1+e_3+e_4}}{(n)_4}\right),
\enas
with
\beas
E_{\alpha,\beta}=\{(e_3,e_4,e_1,e_2): {\bf 1}(e_3=0)+{\bf 1}(e_4+1=0)+{\bf 1}(e_1=0)+{\bf 1}(e_2+1=0)=\alpha\}.
\enas
Letting
\bea \label{def:pk}
p_k(d)= \sum_{b \in \{0,1\}} \frac{(m+1-b)_{k-d}(m+b)_d}{(n)_k}
\ena
we may write
\beas
w_{\alpha,\beta}=\sum_{E_{\alpha,\beta}} p_4(1+e_1+e_2)p_4(1+e_3+e_4).
\enas

Writing out the sets $E_{\alpha,\beta}$ required, we have that
\beas
E_{4,0}=\{(0,1,0,1)\},
\enas
\beas
E_{3,1}=\{(0,1,0,0),(0,1,1,1), (0,0,0,1),(1,1,0,1)\},
\enas
\beas
E_{2,2}=\{(0,0,0,0), (0,0,1,1), (0,1,1,0), (1,0,0,1), (1,1,0,0), (1,1,1,1)\},
\enas
\beas
E_{1,3}=\{(0,0,1,0),(1,0,0,0),(1,0,1,1),(1,1,1,0)\} \qmq{and}
E_{0,4}=\{(1,0,1,0)\}.
\enas
As $p_k(d)=p_k(k-d)$, and since
$E_{\beta,\alpha}$ is obtained from $E_{\alpha,\beta}$ by replacing each entry of each vector by its binary complement,
\beas
w_{\alpha,\beta}=\sum_{E_{\alpha,\beta}}p_4(3-e_1-e_2)p_4(3-e_3-e_4)=\sum_{E_{\beta,\alpha}}p_4(1+e_1+e_2)p_4(1+e_3+e_4)
=w_{\beta,\alpha}.
\enas

Now we calculate the required weights. As the one element of $E_{4,0}$ satisfies $(e_1+e_2,e_3+e_4)=(1,1)$, we have
\beas
w_{4,0}=p_4(2)^2=\left[\frac{2(m)_2(m+1)_2}{(n)_4}\right]^2.
\enas
Next, as the four elements in $E_{3,1}$ yield the vectors $(1,0),(1,2),(0,1),(2,1)$ for
$(e_1+e_2,e_3+e_4)$, we have
\beas
w_{3,1}&=&2p_4(2)p_4(1)+2p_4(2)p_4(3)=4p_4(2)p_4(1)=8\left[\frac{(m)_2(m+1)_2}{(n)_4}\right]
\left[ \frac{(m+1)_3 m+(m)_3(m+1) }{(n)_4} \right].
\enas
Lastly, as the six vectors in  $E_{2,2}$ yield $(0,0),(0,2),(1,1),(1,1),(2,0),(0,0)$ for $(e_1+e_2,e_3+e_4)$,
we have
\beas
w_{2,2}&=&2p_4(1)^2+2p_4(1)p_4(3)+2p_4(2)^2=4p_4(1)^2+2p_4(2)^2\\
&=&4 \left[ \frac{(m+1)_3 m+(m)_3(m+1) }{(n)_4} \right]^2
+ 8 \left[\frac{(m)_2(m+1)_2}{(n)_4}\right]^2.
\enas
Now using $w_{\alpha,\beta}=w_{\beta,\alpha}$,  applying Corollary \ref{cor:spectral.4} in
(\ref{sec.sum.xi-n-4:exp}) we obtain
\begin{eqnarray*}
& &\sum_{a,b \in \{0,1\}} g_{[(1,0), (2, 1); (3, 0), (4, 1)],
{\bf n}(b, a)}^{(m,m+1)} \\
&=& \frac{1}{2}(1+ 6 \lambda_{n,2,{\bf n}_{m,m+1} } + \lambda_{n,4,{\bf n}_{m,m+1} })\left(\frac{(m)_2(m+1)_2}{(n)_4}\right)^2 \\
& & +(1 - \lambda_{n,4,{\bf n}_{m,m+1} })\left(\frac{(m)_2(m+1)_2}{(n)_4}\right)
\left( \frac{(m+1)_3 m+(m)_3(m+1) }{(n)_4} \right) \\
& & + \frac{1}{16}(1-2\lambda_{n,2,{\bf n}_{m,m+1} }+\lambda_{n,4,{\bf n}_{m,m+1} })
\left(4 \left[ \frac{(m+1)_3 m+(m)_3(m+1) }{(n)_4} \right]^2
+ 8 \left[\frac{(m)_2(m+1)_2}{(n)_4}\right]^2 \right).
\end{eqnarray*}
Simplification yields equation (\ref{g:xi-42}).

Next we demonstrate that
\bea  \nn
\sum_{a, b \in \{0, 1\}} \sum_{r \ne t} \left(
g_{[(1, 0), (2, 1); (1, 0), (3, 1)], {\bf
n}(b, a)}^{(r,t)} +g_{[(1, 1), (2, 0); (1, 1), (3, 0)], {\bf n}(b, a)}^{(r,t)} \right.  \\
 + \left.  g_{[(1, 0), (2, 1); (1, 1), (3, 0)],{\bf n}(b, a)}^{(r,t)} + g_{[(1, 1), (2, 0); (1, 0), (3, 1)], {\bf n}(b, a)}^{(r,t)} \right)\nn \\
 = \frac{4(m+1)_2^2(2m-1)^2}{(n)_3^2}.\label{last.4.xi-n-3:exp}
\ena

Starting with the first term, indexed by $[(1, 0), (2, 1); (1, 0), (3, 1)]$, arguing as for the case $p=4$, with $\alpha+\beta=3$ as in (\ref{use.for.3}) we obtain the conditioning event $(X_{01},X_{02},X_{03})=(0, e_2+1,1+e_1)$ and
may write
\bea \label{like.n.4}
\sum_{a,b, \in \{0,1\}} g_{[(1, 0), (2, 1); (1, 0), (3, 1)], {\bf
n}(b, a)}^{(r,t)} =  \sum_{\alpha=0}^3 f_{\alpha,\beta,{\bf n}_{m,m+1}} w_{\alpha,\beta},
\ena
where
\beas
E_{\alpha,\beta}=\{(e_1,e_2): 1+{\bf 1}(e_2+1=0)+{\bf 1}(e_1+1=0)=\alpha\}=\{(e_1,e_2): 1+{\bf 1}(e_2=1)+{\bf 1}(e_1=1)=\alpha\}
\enas
and, with $p_k(d)$ as in (\ref{def:pk}),
\beas
w_{\alpha,\beta}=\sum_{E_{\alpha,\beta}}p_3(1+e_1)p_3(1+e_2).
\enas
Explicitly, the required sets are given by
\beas
E_{3,0}=\{(1,1)\} \quad E_{2,1}=\{(0,1),(1,0)\}, \quad E_{1,2}=\{(0,0)\} \qmq{and} E_{0,3}=\emptyset,
\enas
and calculating the weights, using $p_3(1)=p_3(2)=(m+1)_2(2m-1)/(n)_3$, we obtain
\beas
w_{3,0}= \frac{(m+1)_2^2(2m-1)^2}{(n)_3^2}, \quad
w_{2,1}= \frac{2(m+1)_2^2(2m-1)^2}{(n)_3^2} \qmq{and}
w_{1,2} = \frac{(m+1)_2^2(2m-1)^2}{(n)_3^2}.
\enas
Hence, applying Corollary \ref{cor:spectral.4} in (\ref{like.n.4}) yields
\bea \label{sum3alpha.11}
\lefteqn{\sum_{a,b, \in \{0,1\}} g_{[(1, 0), (2, 1); (1, 0), (3, 1)], {\bf
n}(b, a)}^{(r,t)}}\\
&=&   \frac{(m+1)_2^2(2m-1)^2}{(n)_3^2} \left( \frac{1}{8}(1+3\lambda_{n,1,{\bf n}_{m,m+1}} + 3\lambda_{n,2,{\bf n}_{m,m+1} } + \lambda_{n,3,{\bf n}_{m,m+1}}) \right.\nn \\
&& +\frac{1}{4}(1+\lambda_{n,1,{\bf n}_{m,m+1}} - \lambda_{n,2,{\bf n}_{m,m+1} } - \lambda_{n,3,{\bf n}_{m,m+1}})\nn \\
&&\left. + \frac{1}{8}(1- \lambda_{n,1,{\bf n}_{m,m+1}} - \lambda_{n,2,{\bf n}_{m,m+1} } + \lambda_{n,3,{\bf n}_{m,m+1}})\right). \nn
\ena

Similarly, for the term in (\ref{last.4.xi-n-3:exp}) subscripted by $[(1, 0), (2, 1); (1, 1), (3, 0)]$ the conditioning event
becomes $(X_{01},X_{02},X_{03})=(1, e_2+1,e_1)$, leading to the collection of sets
\beas
E_{3,0}=\emptyset, \quad E_{2,1}=\{(1,0)\}, \quad E_{1,2}=\{(0,0),(1,1)\}, \qmq{and} E_{0,3}=\{(0,1)\},
\enas
that give rise to the same weights as for the terms in (\ref{sum3alpha.11}).
Hence, this term contributes
\bea \label{sum3alpha.12}
&&\frac{(m+1)_2^2(2m-1)^2}{(n)_3^2} \left( \frac{1}{8}(1+\lambda_{n,1,{\bf n}_{m,m+1}} -\lambda_{n,2,{\bf n}_{m,m+1} } -\lambda_{n,3,{\bf n}_{m,m+1}}) \right. \\
&& +\frac{1}{4}(1-\lambda_{n,1,{\bf n}_{m,m+1}} - \lambda_{n,2,{\bf n}_{m,m+1} } +\lambda_{n,3,{\bf n}_{m,m+1}})\nn \\
&&\left. + \frac{1}{8}(1-3\lambda_{n,1,{\bf n}_{m,m+1}} +3 \lambda_{n,2,{\bf n}_{m,m+1} } - \lambda_{n,3,{\bf n}_{m,m+1}})\right), \nn
\ena
and adding to (\ref{sum3alpha.11}) results in
\begin{multline}
\frac{(m+1)_2^2(2m-1)^2}{(n)_3^2} \left( \frac{1}{4}(1+3 \lambda_{n,2,{\bf n}_{m,m+1}})
+\frac{1}{2}(1- \lambda_{n,2,{\bf n}_{m,m+1} }) + \frac{1}{4}(1-\lambda_{n,2,{\bf n}_{m,m+1} }) \right)\\
=\frac{(m+1)_2^2(2m-1)^2}{(n)_3^2}. \label{no.lambda}
\end{multline}

For the term in the sum subscripted by $[(1, 1), (2, 0); (1, 1), (3, 0)]$ the conditioning event
is $(X_{01},X_{02},X_{03})=(0,e_1,e_2)$, resulting in sets
\beas
E_{3,0}=\{(0,0)\}, \quad E_{2,1}=\{(0,1),(1,0)\}, \quad E_{1,2}=\{(1,1)\} \qmq{and} E_{0,3}=\emptyset
\enas
and a term, therefore, agreeing with (\ref{sum3alpha.11}). Lastly, for the term subscripted by $[(1, 1), (2, 0); (1, 0), (3, 1)]$ the conditioning event is $(X_{01},X_{02},X_{03})=(1, e_1,1+e_2)$,
the collection of sets becomes,
\beas
E_{3,0}=\emptyset, \quad E_{2,1}=\{(0,1)\}, \quad E_{1,2}=\{(0,0),(1,1)\}, \qmq{and} E_{0,3}=\{(1,0)\},
\enas
and we thus obtain a term agreeing with (\ref{sum3alpha.12}). Hence these last two terms
also contribute the factor (\ref{no.lambda}), and summing over the cases
$(r,t)=(m,m+1)$ and $(r,t)=(m+1,m)$ yields (\ref{last.4.xi-n-3:exp}).

Lastly, we show that for $r,t \in \{m,m+1\},r \not = t$,
\bea \label{g.n221}
\sum_{a,b \in \{0,1\}} g_{[(1,0),(2,1);(1,0),(2,1)], {\bf n}(b, a)}^{(r,t)} =\left(1+2\lambda_{n,1,{\bf n}_{m,m+1}} + \lambda_{n,2,{\bf n}_{m,m+1} } \right)\frac{(m+1)_2^2}{(n)_2^2}.
\ena
In this case the conditioning event is simply $(X_{01},X_{02})=(0,0)$,
and selecting one on and one off switch in each of the stages $r$ and $t$ yields the weight $w_{2,0}=p_2(1)^2=4(m+1)_2^2/(n)_2^2$; applying Corollary \ref{cor:spectral.4} now completes the verification of (\ref{g.n221}).

\section*{Bibliography}
\begin{enumerate}

\item \label{BRS} Baldi, P., Rinott, Y.  and Stein, C.  (1989).
A normal approximation for  the number of local maxima of a random
function on a graph, In {\em Probability, Statistics and Mathematics,
Papers in Honor of Samuel Karlin.}  T. W. Anderson, K.B. Athreya
and D. L. Iglehart eds., Academic Press , 59-81.


\item \label{GhGo} Ghosh, S. and Goldstein, L. (2009). Concentration of measures via size biased couplings. {\em Probability Theory and Related Fields, to appear}

\item \label{Chen} Chen, L.H.Y., Goldstein, L., and Shao, Q.M. (2010). Normal Approximation by Stein's Method.
Springer.

\item \label{Gold}
Goldstein, L. (2005).
Berry-Esseen bounds for combinatorial central limit theorems
 and pattern occurrences, using zero and size biasing.
{\em J. Appl. Probab.} {\bf 42}, 661--683.

\item \label{GR}
Goldstein, L. and Rinott, Y. (1996). On multivariate normal
approximations by Stein's method and size bias couplings. {\em J.
Appl. Prob.} {\bf 33}, 1-17.

\item \label{GolPen}
Goldstein, L. and Penrose, M. (2009). Normal approximation for coverage
 models over binomial point processes. {\em Ann. Appl. Prob}. {\bf
 20}, 696 - 721.


\item \label{rrz} Rao, C., Rao, M., and Zhang, H. (2007).
One Bulb? Two Bulbs? How many bulbs light up? A discrete probability
problem involving dermal patches. {\em Sanky$\bar{a}$}, {\bf 69},
137-161.

\item \label{ShSu} Shao, Q. M. and Su, Z. (2005). The Berry-Esseen bound for character ratios.
{\em Proceedings of the American Mathematical Society}, {\bf 134}, 2153-2159.

\item \label{Steinchar} Stein, C. (1981) Estimation of the mean of a multivariate normal distribution.  {\em Ann. Statist.} {\bf 9}, 1135-1151.

\item \label{Stein86} Stein, C. (1986). Approximate Computation of Expectations. IMS, Hayward, California.

\item \label{zl} Zhou, H. and Lange, K. (2009). Composition Markov chains of multinomial type. {\em Advances in Applied Probability}, {\bf 41}, 270 - 291.

 \end{enumerate}

\end{document}